\documentstyle[twoside,12pt]{amsart}

\def\myemail{phung@@msri.org}
\def\myperemail{{phung@@ioit.ncst.ac.vn}}

\def\myaddress{Mathematical Sciences Research Institute, Berkeley, CA 94720.}
\def\myperaddress{Hanoi Institute of Mathematics, P.O.Box 631, 10000 Boho, Hanoi}
\def\mythanks{The author should like to thank Professor C. Procesi for explaining him invariant theory. He also should like to thank Professor K. Goodearl for useful discussions. 

The first part of the work was done at the Max-Planck Intitut f\"ur Mathematik, Bonn and appeared as Preprint MPI-99/12. The work was completed during the author's stay at the Mathematical Sciences Research Institute, Berkeley. The author would like to thank these Institutes for the excellent working condition and financial support.}
\def\myabstract{For a Hecke operator $R$, one defines the matrix bialgebra $\E_R$, which is considered as the function algebra on the quantum space of endomorphisms of the quantum space associated to $R$. One generalizes this notion, defining the function algebra $\M_{RS}$ on the quantum space of homomorphisms of two quantum spaces associated to two Hecke operators $R$ and $S$ respectively. $\M_{RS}$ can be considered as a quantum analogue (or a deformation) of the function algebra on the variety of matrices of a certain degree. We provide two realiztions of $\M_{RS}$ as a quotient algebra and as a subalgebra of a tensor algebra, whence derive interesting informations about $\M_{RS}$, for instance the Koszul property, a formula for computing the Poincar\'e series. On $\M_{RS}$ coact the bialgebras $\E_R$ and $\E_S$. We study the two-sided ideals in $\M_{RS}$, invariant with respect to these actions, in particular, the determinantal ideals. We prove analogies of the fundamental theorems on invariant theory for these quantum groups and quantum hom-spaces.}

\def\SS{{\cal S}}
\def\S{{\sf S}}
\def\F{{\sf F}}
\def\Ff{{\frak F}}
\def\Gg{{\frak G}}

\def\myauthor{Ph\`ung H{{\accent"5E o}\kern-.28em\raise.2ex\hbox{\char'22}\kern-.20em} H{a\kern-.370em\raise.16ex\hbox{\char'47}\kern.1em}i}
\def\myAUTHOR{ PH\`UNG H{{\accent"5E O}\kern-.38em\raise.8ex\hbox{\char'22}\kern-.12em}  H{A\kern-.46em\raise.80ex\hbox{\char'47}\kern.18em}I}

\def\amshead{
\title[Quantum hom-spaces, invariant theory and determinantal ideals]{Realizations of quantum hom-spaces, invariant theory and quantum determinantal ideals}
\author{ PH\`UNG H{{\accent"5E O}\kern-.38em\raise.8ex\hbox{\char'22}\kern-.12em}  H{A\kern-.46em\raise.80ex\hbox{\char'47}\kern.18em}i}
\address{\myperaddress}
\curraddr{\myaddress}
\subjclass{Primary 16W30,17B37 , Secondary 17A45, 17A70}
\keywords{Hecke operators, quantum hom-space, quantum determinant, invariant theory}
\email{\myperemail\ and \myemail}
\begin{abstract}\myabstract\end{abstract}
\maketitle }

\def\G{{\sf G}}
\def\lora{\longrightarrow}
\def\ot{\otimes}

\def\loma{\longmapsto}
\def\Vn{{V^{\ot n}}}
\def\Vns{V^{*\ot n}}

\def\Vsn{{V^{*\ot n}}}

\def\si{\sigma}
\newcommand{\bbas}{\begin{eqnarray*}}
\newcommand{\eeas}{\end{eqnarray*}}
\newcommand{\bbar}{\begin{array}}
\newcommand{\eear}{\end{array}}
\newcommand{\bbs}{\begin{displaymath}}
\newcommand{\ees}{\end{displaymath}}
\newcommand{\bb}{\begin{equation}}
\newcommand{\eqbb}{\begin{equation}}
\def\ee{\end{equation}}
\def\eqee{\end{equation}}
\def\eea{\end{eqnarray}}

\def\bba{\begin{eqnarray}}
\newtheorem{thm}{Theorem}[section]
\newtheorem{lem}[thm]{Lemma}

\newtheorem{cor}[thm]{Corollary}
\newtheorem{pro}[thm]{Proposition}

\def\Ker{\mbox{\sf Ker\,}}
\def\Im{\mbox{\sf Im\,}}

\def\Mor{\mbox{\sf Mor\,}}
\def\Hom{\mbox{\sf Hom\,}}

\def\End{\mbox{\sf End\,}}

\def\H{{\cal H}}
\def\Hb{{\sf H}}
\def\O{{\cal O}}

\def\A{{\sf A}}
\def\AA{{\cal A}}

\def\N{{\sf N}}

\def\T{{\sf T}}
\def\V{{\cal V}}
\def\P{{\cal P}}
\def\D{{\cal D}}
\def\E{{\sf E}}

\def\eee{\rule{.75ex}{1.5ex}\\[1ex]}

\def\proof{{\it Proof.\ }}
\newcommand{\va}{\varepsilon}

\newcommand{\Nat}{\mbox{\sf Nat\,}}

\renewcommand{\dim}{\mbox{\sf dim\,}}

\def\rref#1{(\ref{#1})}
\def\M{{\sf M}}

\font\Fraktur=eufm10 scaled\magstep1          
   \newcommand{\fraktur}[1]{\mbox{\Fraktur #1}}  %
   \font\Fraktu=eufm7 scaled\magstep1            
   \newcommand{\fraktu}[1]{\mbox{\Fraktu #1}}    %
   \font\Frakt=eufm5 scaled\magstep1             
  \newcommand{\frakt}[1]{\mbox{\Frakt #1}}      %
   \def\frak#1{\mathchoice{\fraktur {#1}}            
                        {\fraktur {#1}}            
                        {\fraktu {#1}}             
                        {\frakt {#1}}  }           

\newcommand{\Ss}{\frak S}

\font\Bbbm=msbm10 scaled\magstep1
  \def\BBbm#1{\mbox{\Bbbm #1}}
\font\Bbbb=msbm7 scaled\magstep1
  \def\BBbb#1{\mbox{\Bbbb #1}}
\font\Bbbn=msbm5 scaled\magstep1
  \def\BBbn#1{\mbox{\Bbbn #1}}

\def\Bbb#1{\mathchoice{\BBbm {#1}}
                    {\BBbm {#1}}
                    {\BBbb {#1}}
                    {\BBbn {#1}} }
\newcommand{\bK}{{\Bbb K}}

\newcommand{\bN}{{\Bbb N}}

\def\db{{\mathchoice{\mbox{\sf db}}
                    {\mbox{\sf db}}
                    {\mbox{\scriptsize\sf db}}
                    {\mbox{\tiny\sf db}} }}

\def\ev{{\mathchoice{\mbox{\sf ev}}
                    {\mbox{\sf ev}}
                    {\mbox{\scriptsize\sf ev}}
                    {\mbox{\tiny\sf ev}} }}

\def\id{{\mathchoice{\mbox{\sf id}}
                    {\mbox{\sf id}}
                    {\mbox{\scriptsize\sf id}}
                    {\mbox{\tiny\sf id}} }}

\def\op{{\mathchoice{\mbox{\rm op}}
                    {\mbox{\rm op}}
                    {\mbox{\scriptsize\rm op}}
                    {\mbox{\tiny\rm op}} }}

\def\part{\vdash}

\def\lam{{\lambda}}

\def\cdotot{} 

\def\Lambda{{\mathchoice{\mbox{\raisebox{.25ex}{$\textstyle\bigwedge$}}}
{\mbox{\raisebox{.25ex}{$\textstyle\bigwedge$}}}
{\mbox{\raisebox{.17ex}{\scriptsize$\textstyle\wedge$}}}
{\mbox{\raisebox{.1ex}{\tiny$\textstyle\wedge$}}}
} }

\def\Rt{{}^t\!R}
\def\St{{}^t\!S}

\begin{document}
\bibliographystyle{plain}
\amshead
\setcounter{section}{-1}
\section*{Introduction}
Let $V$ be a vector space of finite dimension over a field of characteristic zero. There are two ways of interpreting the symmetric tensor algebra $\S(V)$ over $V$. The first one is to consider it as a factor algebra of the tensor algebra $\T(V)$ over $V$. Thus we have a {\it relation} between two tensors, like $a\ot b=b\ot a$. The second way is to consider $\S(V)$ as a {\it subspace} of $\T(V)$ consisting of symmetrized tensors, e.g., for $a,b\in V$, the element $a\ot b+b\ot a$ belongs to the symmetric tensor algebra on $V$. We are thus speaking of two realizations of the symmetric tensor algebras over $V$. The analogous procedure applies also for the exterior (anti-symmetric) tensor algebra.

A Hecke operator $R$ on a vector space $V$ of finite dimension is an invertible operator on $V\ot V$ that satisfies the Yang-Baxter equation and the Hecke equation $(x+1)(x-q)=0$. To a Hecke operator there is associated a quantum space, given in terms of a pair of quadratic algebras (Section \ref{sec2}). The latter algebras are quantum analogues (or deformations) of the symmetric and anti-symmetric tensor algebras over a vector space.

For a quantum space associated to a Hecke operator, the two realizations for its (quantum) symmetric and anti-symmetric tensor algebras was first obtained by Gurevich \cite{gur1}. While the first realization was taken as the definition, the second realizations followed from the general theory of Hecke algebras. In this paper we give the second realization for the (anti-) symmetric tensor algebras on the quantum semi-group of endomorphism associated to a Hecke operator and the quantum hom-space associated to a pair of Hecke operators.

The main idea in giving the second realization is to construct a projector on the tensor power of $V$. For the classical case, it is the (anti-) symmetrizer operator, e.g., $a\ot b\lora (a\ot b+b\ot a)/2$. For the quantum space, it is the quantum (anti-) symmetrizer, constructed in terms of the trivial and signature representations of the Hecke algebras. In Section \ref{sec2}, we recall the definition of the quantum exterior and the quantum symmetric algebras associated to a Hecke operator and their second realization, due to Gurevich (Equation \rref{eq6}). Together, these algebras determine a quantum space. Then, we recall the definition of the matrix bialgebra associated to the Hecke operator, which is considered as the function algebra of the quantum semi-group of endomorphism of the quantum space. Unlike the classical case, where a matrix can also be considered as a vector, the matrix bialgebra generally cannot be defined as a quantum symmetric algebra associated to a Hecke operator. In fact, it can still be defined analogously in terms of a Yang-Baxter operator, but this operator has the minimal polynomial of degree 3. As a result, only an analogue of the quantum anti-symmetrizer was defined (the operator $\Phi^n$ in Equation \rref{eq10}), which is no more a projector. There is not strightforward analogue of the quantum symmetrizer.

There is a simple solution by the following remark. Since the Yang-Baxter operator defining the matrix bialgebra is a tensor product of the ordinary Hecke operator with the inverse of its dual, our operator $\Phi^n$ is a homomorphic image of a Casimir element in $\H_n\ot\H_n$, where $\H_n$ is the Hecke algebra of type $A_n$. Using various dual bases in $\H_n$, we found eigenvalues of $\Phi^n$. So that we can modify it to obtain a projector $\bar\Phi^n$ (Subsection \ref{sect-phi}). Now, the choice of a quantum symmetrizer becomes clear. The remark above also suggest us define an operator $\Psi^n$, which play the role of the quantum symmetrizer. We find its eigenvalues and modify it to get a projector $\bar\Psi^n$. It is the operators $\bar\Psi^n$, $n=1,2,\ldots$, which give the second realization for the matrix bialgebra.

The results can be generalized for the function algebra on the space of homomorphisms of two quantum spaces or quantum hom-space. This quadratic algebra is introduced in Section \ref{sec4}, it is defined in terms of a pair of Hecke operators. Thus, on it coact the matrix bialgebras associated to these Hecke operators. Our second realization for quantum spaces of homomorphisms  implies a quantum analogue of Cauchy's decomposition \rref{eq20}. Although in the classical case, the operator $\bar\Psi^n$ reduces to the ordinary symmetrizer operator, its relationship with Cauchy's decomposition is new. Moreover, the second realization also implies interesting results in invariant theory.

Before describing the invariant theory for quantum groups of type $A$, let me briefly recall the classical theory. Let $M(m,n)$ denote the space of $m\times n$-matrices.
Let $\mu$ denote the matrix multiplication map $\mu:M(m,t)\times M(t,n)\lora M(m,n)$, $\mu(A,B)=AB$, $A\in M(m,t), B\in M(t,n)$. On the variety $M(m,t)\times M(t,n)$ acts the general linear group $GL(t)$, $g(A,B)=(Ag^{-1},gB)$. It this easy to see that elements of an orbit of $GL(t)$ have the same image under $\mu$. The above action of $GL(t)$ induces an action on the polynomial ring on $M(m,t)\times M(t,n)$, $\O(M(m,t)\times M(t,n))$. The classical invariant theory for general linear groups studies the subring of invariant polynomials in $\O(M(m,t)\times M(t,n))$. Let $m^i_j$ the $(i,j)$ coordinate function on $\O(M(m,n))$ and $a^i_j$ be the composition of $\mu$ with $m^i_j$, which are then polynomial functions on $\O(M(m,t)\times M(t,n))$. The first fundamental theorem of invariant theory states that any invariant polynomial on $\O(M(m,t)\times M(t,n))$ can be represented as a polynomial functions on the functions $a^i_j$. The second fundamental theorem states that the relations between the functions $a^i_j$ are exactly the minors of ranks $t+1$ in the matrix $(a^i_j)$. 

Using the associated homomorphism of algebras $\mu^*:\O(M(m,n))\lora \O(M(m,t)\times M(t,n))$, we can reformulate the above theorems as follows:
\begin{enumerate}\item A polynomial in $\O(M(m,t)\times M(t,n))$, invariant under the action of $GL(t)$, is contained in the image of $\mu^*$.
\item The kernel of $\mu^*$ is the ideal in $\O(M(m,n))$, generated by minors of degree $t+1$ in the matrix $(e^i_j)$.\end{enumerate}

The characteristic free proof of these theorems, due to DeConcini-Procesi \cite{cp1}, uses the notion of standard basis and has close relationship with combinatorics. A generalization of these results for standard quantum general linear groups was obtained by Goodearl et.al \cite{gl,glr}. Their proof closely follows DeConcini-Procesi's proof.

In the quantum setting, the variety $M(m,n)$ is replaced by a quantum hom-space. Given two Hecke operators $R,S$, the quantum hom-space associated to $R,S$ is denoted by $\M_{SR}$. On the algebra $\M_{SR}$ coact the bialgebras $\E_R$ on the right and $\E_S$ on the left. In Section \ref{sec5}, we study two sided ideals in $\M_{SR}$ which are invariant with respect to these actions. In the classical case, these ideals were studied by DeConcini-Eisenbud-Procesi \cite{cep1}. We show that there is a one-one correspondence between invariant ideals and diagram ideals in the sense of \cite{cep1} (note that our notation here slightly differes from the notion in \cite{cep1}, a partition is replaced by its conjugate partition). In the classical case, those invariant ideals that correspond to diagram ideals of the form $\langle (1^k)\rangle=\{\lambda|\lam_1\geq k\}$ are determinantal ideals, i.e., generated by minors of degree $k$. In our general case, the notions of quantum determinant and quantum minor is not defined. We show, however, that in the case of standard quantum general linear groups, where these notions are defined, the quantum determinantal ideals, introduced by Goodearl et.al., are precisely those corresponding to $\langle(1^k)\rangle$, for some $k$.

The setting for invariant theory of quantum groups of type $A$ involves three Hecke operators $R,S$ and $T$, of which $R$ is also a Hecke symmetries. The morphism $\mu^*$ mentioned above becomes an algebra morphism
\bbs \mu^*:\lora \M_{TS}\lora \M_{TR}\ot \M_{RS}.\ees
The corresponding quantum group is the Hopf algebra associated to $R$, $\Hb_R$. This Hopf algebra coacts on the source and the target of $\mu^*$ and the formulation of the first and the second fundamental theorems can be made analogously as in the classical case. Since $\Hb_R$ is not commutative, there are more than one coactions of it on $\M_{RS}$, which yield different versions of the fundamental theorems.

The method of our proof is new. It relies mainly on the second realization of function algebras on quantum spaces homomorphisms. In the case of standard quantum general linear groups, our result is precisely those obtained in \cite{gl,glr}. On the other hand, our assumption about the quantum groups also covers the case of standard quantum general linear supergroups. Thus, we have particularly proved the fundamental theorems for quantum general linear supergroups, which have as a special case the fundamental theorems for general linear supergroups. Unlike the case of (quantum) general linear groups, in the super case, the kernel of $\mu^*$ is not generated by (quantum) minors. It is an interesting problem to study such ideals.

\section{Preliminaries}\label{sec1}
Throughout this paper, we work over an algebraically closed field $\bK$ of characteristic zero.
\subsection{Partitions.}\label{sec1.1} A partitions  $\lam$ of $n\in\bN$ is a sequence $\lam=(\lam_1,\lam_2,...)$ of non-increasing non-negative integers, whose sum is $n$, we write $\lam\part n$ or $|\lam|=n$. The maximal number $r$, for which $\lam_r\neq 0$, is called the length of $\lam$. The diagram $[\lam]$ associated to $\lam$ is a matrix, whose first row contains $\lam_1$ elements, called nodes, second row contains $\lam_2$ elements, and so on. The conjugated to $\lam$, denoted by $\lam'$, is the one, whose diagram $[\lam']$ is obtained from $[\lam]$ by rotating it $180^\circ$ along its diagonal.

A standard $\lam$-tableau is the diagram $[\lam]$, filled by numbers $1,2,...,|\lam|$, in such a way that they increase along rows and columns. The number of standard tableaux is denoted by $d_\lam$.

\subsection{Symmetric Groups.}\label{sec1.2} The symmetric group $\Ss_n$ consists of permutations of the set $1,2,\ldots,n$. It can be regarded as the group generated by transposition $v_i=(i,i+1)$, $1\leq i\leq n-1$, subject to the relations: $v_i^2=1, v_iv_{i+1}v_i=v_{i+1}v_iv_{i+1}$ and $v_iv_j=v_jv_i$ if $|i-j|\geq 2$. The length $l(w)$ of an element $w$ is the minimal length of the words in $v_i$ expressing $w$. It equals the number of pairs $1\leq i<j\leq n$ for which $iw>jw$.

\subsection{Hecke Algebras \cite{dj1}.}\label{sec1.3} The Hecke algebra $\H_n=\H_{n,q}$ is a $q$-analogue of the group algebra $\bK[\Ss_n]$. It is generated over $\bK$ by 1 and the elements $T_i,1\leq i\leq n-1$, subject to the relations $T_i^2=(q-1)T_i+q$, $T_iT_{i+1}T_i=T_{i+1}T_iT_{i+1}$ and $T_iT_j=T_jT_i$ if $|i-j|\geq 2$. When $q=1$, $\H_{n,1}$ reduces to $\bK[\Ss_n]$, $T_i\lora v_i$. $\H_n$ has a basis consisting of $T_w, w\in\Ss_n$, $T_1:=1$, $T_{v_i}:=T_i$, $T_wT_u=T_{wu}$ if $l(w)+l(u)=l(wu)$.

We shall always assume that $q^n\neq 1, \forall n>1$. In this case $\H_{n,q}$ is semisimple.

The embedding $\H_l\ot \H_m\lora \H_{l+m}$, mapping $\H_l\ni T_i$ to $T_i\in\H_{l+m}$ and $\H_m\ni T_j$ to $I_{m+j}\in\H_{l+m}$, is called the standard embedding.

\subsection{Representations of the Hecke Algebras \cite{dj1}.} \label{sec1.4}Representations of $\H_{n,q}$, for $q$ not being root of unity, can be parameterized by partitions of $n$. Let $S_\lam$ be the simple representation of $\H_n$, corresponding to $\lam$, then the dimension of $S_\lam$ is $d_\lam$ ($d_\lam$ is defined in \ref{sec1.1}). Let $\AA_\lam$ be the block, i.e. a minimal two sided ideal, in $\H_n$, that corresponds to $\lam\part n$. Let $E_\lam^{ij}, 1\leq i,j\leq d_\lam$ be a basis of $\AA_\lam$ such that 
$E_\lam^{ij}E_\lam^{kl}=\delta^j_kE_\lam^{il}.$ Thus, $E_\lam^{ii}$ are mutually orthogonal primitive idempontents of $\H_n$. Notice that for $\lam\neq \mu$, $E_\lam^{ij}E_\mu^{kl}=0$.

We shall also consider the algebra $\H_n^{\rm op}$. Its simple comodules are canonically identified with ${S_\lam}^*$, the dual vector space to $S_\lam$: if $\phi_\lam$ is the representation of $\H_n$ on $S_\lam$, then the representation of $\H_n^{\rm op}$ on ${S_\lam}$ is given by $\bar\phi_\lam(w):=\phi_\lam(w)^*.$

\subsection{A bilinear form on $\H_n$ \cite{dj1}.}\label{sec1.5} There exists a non-degenerate, symmetric, associative bilinear form on $\H_n$, defined as follows: $(T_u,T_w):=q^{l(u)}\delta^{u^{-1}}_w.$ Thus, we see that $\{q^{-l(w)}T_{w^{-1}}, w\in\Ss_n\}$ is the dual basis to $\{T_w, w\in\Ss_n\}$ with respect to this bilinear form. Therefore, the Casimir element $\sum_{w\in\Ss_n}q^{-l(w)}T_w\ot T_{w^{-1}}$ is central:
\bbs \sum_{w\in\Ss_n}q^{-l(w)}T_iT_w\ot T_{w^{-1}}=\sum_{w\in\Ss_n}q^{-l(w)}T_w\ot T_{w^{-1}}T_i.\ees

Since $\AA_\lam$ is simple, the bilinear form restricted on $\AA_\lam$ should satisfy $(E_\lam^{ij},E_\lam^{kl})=\delta^j_k\delta^i_l k_\lam$, for certain coefficient $k_\lam\neq 0$. 
Thus $\{E_\lam^{ij},\lam\part n, 1\leq i,j\leq d_\lam\}$ and $\{k_\lam^{-1}E_\lam^{ji},\lam\part n, 1\leq i,j\leq d_\lam\}$ are dual bases. Therefore
\bba\label{eq0} 
\sum_{w\in\Ss_n}q^{-l(w)}T_w\ot T_{w^{-1}}=\sum_{\lam\part n\atop 1\leq i,j\leq d_\lam}k^{-1}_\lam E_\lam^{ij}\ot E_\lam^{ji}.\eea
Denote $\Pi_\lam:=\sum_{ 1\leq i,j\leq d_\lam}k^{-1}_\lam E_\lam^{ij}\ot E_\lam^{ji}.$ Then
\bbs \Pi_\lam^2=d_\lam k_\lam^{-1}\Pi_\lam.\ees
The coefficients $k_\lam$ can be explicitly given (cf. \cite{ph98})
\bbs k_\lam=q^{\sum_i{\lambda_i}(i-1)}\prod_{k=1}^n\frac{1}{[r_\lambda(k)+r]_q}\prod\frac{[\lambda_i-\lambda_j+j-i]_q}{[j-i]_q}.\ees
\subsection{Hecke Operators.}\label{sec1.6} Let $V$ be a finite dimensional vector space over a field $\bK$ and $R$ be a invertible operator on $V\ot V$. $R$ is called Hecke operator if the following conditions are satisfied:
\bbas && R_1R_2R_1=R_2R_1R_2\quad \mbox{where } R_1:=R\ot \id_V, R_2:=\id_V\ot R,\\
&& (R+1)(R-q)=0.\eeas

A Hecke operator induces a right representation $\rho$ of $\H_{n,q}$ on $\Vn$, for $n>1$: $\rho(T_i)=R_i:=\id_V^{\ot i-1}\ot R\ot\id_V^{\ot n-i-1}.$ For any $w\in\Ss_n$, we denote $R_w:=\rho(T_w)$. 

\subsection{Hopf algebras.}\label{sec1.7} We assume that the reader is familiar with the notions of bialgebras and Hopf algebras and their (co)modules. The reader may consult \cite{sweedler1} or \cite[Chapter1]{manin1} for basic notions of bialgebras and Hopf algebras. For a coalgebra $C$, $C^{\rm cop}$ denotes $C$ with the opposite coproduct. Similary, for a bialgebra $B$, $B^{\rm cop}$ denotes the $B$ with the opposite product and coproduct, it is a bialgebra too. If $M$ is a right $C$ comodule then $M^*$ is a left $C$-comodule in a canonical way, namely, if $\delta$ denotes the coaction of $C$ on $M$, then the left coaction $\lam$ of $C$ on $M^*$ is given by $\lam(\phi)(m)=\phi(\delta(m))$ for all $m\in M,\phi\in M^*$. Hence, $M^*$ is a right $C^{\rm cop}$-comodule. 

\subsection{The dual space to a tensor product.}\label{sec1.8} To a vector space $V$, the dual vector space is defined to be $V^*:=\Hom(V,\bK)$. If $V$ is finite dimensional over $\bK$, then so is $V^*$ and they have the same dimension. For finite dimension vector spaces, there is the following equivalent definition of dual vector spaces which is more suitable for further generalization (cf. \cite{deligne82}). The dual vector space to a vector space $V$ (of finite dimension) is a pair $(V^*,\ev_V:V^*\ot V\lora \bK)$, such that there exists a linear map $\db_V:\bK\lora V\ot V^*$, satisfying the following conditions ($\id_V$ is the identity map on $V$):
\bbs (\ev\ot \id_{M^*})(\id_{M^*}\ot \db)=\id_{M^*},\quad (\id_M\ot\ev)(\db\ot \id_M)=\id_M.\ees
The dual vector space is determined uniquely up to an isomorphism by these conditions. For the generalization for monoidal categories see \cite[Chap.~1]{deligne82} or \cite{kassel}. Notice that $\ev$ plays the role of the pairing between $V$ and $V^*$.

Although the vector space $(V\ot W)^*:=\Hom(V\ot W,\bK)$ is uniquely determined, there are more than one way of specifying a basis for this space based on given bases on $V$ and $W$. This can be seen from the point of view of the second definition above as the existence of more than one choices of the $\bK$-linear mapping $\ev_{V\ot W}$. More precisely, one can set $(V\ot W)^*=V^*\ot W^*$, which is the usual way, or $(V\ot W)^*=W^*\ot V^*$ which is actually a more standard way. In this paper we shall use both ways of identifying $(V\ot W)^*$. Note that for longer tensor products, we shall identify their duals in the above two ways (although there are more), namely $(V\ot W\ot\cdots\ot U)^*=V^*\ot W^*\ot\cdots\ot U^*$ or $(V\ot W\ot\cdots\ot U)^*=U^*\ot\cdots\ot W^*\ot V^*$. For the case of $V^{\ot n}=V\ot V\ot\cdots \ot V$, the above notation may be confusing, so we shall use the notations $V^{*\ot n}$ and $V^{\ot n*}$ to refer to the first and the second identification, respectively.

One of the reasons to specify the dual spaces is for describing the matrices of adjoint operators. Recall that to each operator $R:U\lora V$ on finite dimensional vector spaces, there corresponds an adjoint operator $R^*:U^*\lora V^*$, which is uniquely determined. For operators on tensor products, the determination of an adjoint operator obviously depends on the choice of the specification of dual spaces. In this paper we shall use the following notation for the matrices of adjoint operators: $\Rt$ for the usual identification $(V\ot W\ot\cdots\ot U)^*=V^*\ot W^*\ot\cdots\ot U^*$ and $R^*$ for the ``standard'' identification $(V\ot W\ot\cdots\ot U)^*=U^*\ot\cdots\ot W^*\ot V^*$.

\section{Matrix Quantum Semigroups of Type $A$}\label{sec2}
Let $\bK$ be an algebraically closed field of characteristic zero, which will be fixed through out this paper. Let $q\in\bK^\times$ which is not a root of unity of order greater that 1. $q$ will also be fixed through out the paper.

\subsection{Quadratic Algebras}\label{sec2.1}

Let $V$ be a vector space over $\bK$ of finite dimension. Let $R$ be a subspace of $V\ot V$. Let $\A=\A(V,R)$ be the quotient algebra of the tensor algebra on $V$ by the two-sided ideal generated by elements of $R$:
 \bbs \A(V,R):=\T(V)/(R).\ees
Such an algebra is called quadratic algebra \cite[Chap.~1]{manin1}. $R$ is called the space of relations of $\A$. The two-sided ideal, generated by $R$, is  usually denoted by $R(\A)$. Since the relations on $\A$ are homogeneous, $\A$ inherits a grading from $\T(V)$
\bbas &&\A=\bigoplus_{n=0}^\infty \A_n, \quad \A_n=V^{\ot n}/R^n(\A)\\
 &&     R^n(\A)=\sum_{i=1}^{n-1}R(\A)^n_i,\quad R(\A)^n_i:=V^{i-1}\ot R\ot V^{n-i-1},\ 1\leq i\leq n-1.\eeas

The Poincar\'e (or Hilbert) series of $\A$ is by definition the formal power series $P_\A(t):=\sum_{k=0}^\infty \dim_\bK \A_k.$ The dual quadratic algebra to $\A$, $\A^!$ is defined to be $\A^!:=\T(V^*)/(R(\A)^\perp).$ It is easy to see that $(\A^!_n)^*=\bigcap_{i=1}^{n-1}R(\A)^n_i.$

\vskip1ex

\noindent{\it Example.} Assume that $V$ has dimension $d$. Then $\T(V)$ is canonically isomorphic to the free non-commutative algebra with $d$ generators: $k\langle x_1,x_2,\ldots,x_d\rangle $. The polynomial ring in $d$ indeterminates is the quotient of this algebra by the ideal generated by elements of the form $x_ix_j-x_jx_i$, thus a quadratic algebra. The Poincar\'e series of the polynomial ring is equal to $(1-t)^{-d}$ as a formal series. The dual quadratic algebra is the polynomial algebra on anti-commuting indeterminates $\Lambda(\xi^1,\xi^2,\ldots,\xi^d)$, it is the quotient of the free non-commutative algebra $k\langle\xi^1,\xi^2,\ldots,\xi^d\rangle $ by the ideal generated by elements of the form $\xi^i\xi^j+\xi^j\xi^i$.

\subsection{Quantum Spaces and Quantum Endomorphism Rings}\label{sec2.2}
Let $R=R_q$ be a Hecke operator on $V$, where $q\in\bK^\times$ will be assumed not to be a root of unity of order greater than 1. We define the following quadratic algebras
\bba &&\S=\S_R:=\T(V)/(\Im(R-q)),\label{eq1}\\
&& \Lambda=\Lambda_R:=\T(V)/(\Im(R+1)).\label{eq2}\eea
 $\S_R$ and $\Lambda_R$ are considered as the function algebra and the exterior algebra on  a quantum space. The function algebra on the quantum semi-group of endomorphisms of this quantum space is defined to be a quadratic algebra on $V^*\ot V$, given by
\bba\label{eq3}\E=\E_R:=\T(V^*\ot V)/(\Im(\bar R-1)),\eea
where $\bar R:=s_{(23)}(R^{*-1}\ot R)s_{(23)}$, acting on $(V^*\ot V)^{\ot 2}$, $s_{(23)}$ interchanges the second and the third components in a tensor product, (see \cite[Section 1.2]{ph97}). One can check that 
\bbs R(\E)=s_{(23)}\left(R(\S)^\perp\ot R(\S)\oplus R(\Lambda)^\perp\ot R(\Lambda)\right).\ees
Using $s_{(23)}$ we shall identify $\bar R$ with $\Rt^{-1}\ot R$ acting on $V^{*\ot 2}\ot V^{\ot 2}$.

Let $\ev_V$ denote the linear map $\ev_V:V^*\ot V\lora k$, $\phi\ot x\loma \phi(x)$. Then there exists uniquely a morphism $\db_V:\bK\lora V\ot V^*$, subject to the following conditions:
\bbs (\ev\ot \id_{V^*})(\id_{V^*}\ot \db)=\id_{V^*},\quad (\id_V\ot\ev)(\db\ot \id_V)=\id_V.\ees
In fact, if $x_1,x_2,\ldots,x_n$ form a basis of $V$ and $\xi^1,\xi^2,\ldots,\xi^n$ form a basis of $V^*$, such that $\ev_V(\xi^i\ot x_j)=\delta_j^i$, then $\db_V$ is given by $\db_V(1)=\sum_k\xi^k\ot x_k$. Conversely, the dual vector space $V^*$ can be given in terms of the linear map $\ev_V$ and $\db_V$. Then $\T_n\cong \Vns\ot \Vn$ is a coalgebra with the coproduct 
\bbs \Delta_n=\id\ot\db_{\Vn}\ot\id:\Vns\ot \Vn\lora \Vns\ot \Vn\ot \Vns\ot \Vn.\ees
The direct sum of $\Delta_n$ defines a coalgebra structure on $\T$ making it a bialgebra. It turns out that the ideal generated by $R(\E)$ is a biideal (cf. \cite[Chap.~2]{manin1}), hence $\E=\T/R(\E)$ is a bialgebra too. Set $e^i_j:=\xi^i\ot x_j$. Then $\{e^i_j:=\xi^i\ot x_j, 1\leq i,j\leq d\}$ is a basis of $V^*\ot V$. Then the coproduct on $\E$ is given by
\bbs \Delta(e^i_j)=\sum_ke^i_k\ot e^k_j.\ees
We shall also write $\Delta(E)=E\dot\ot E,$ for $E=(e^i_j)$.

The relation on $\E$ can be written in terms of the matrix $E$ as follows (cf \cite{ph97}):
\bba\label{eq3.1} RE_1E_2=E_1E_2R, \quad\mbox{where $E_1:=E\ot \id(d), E_2:=\id(d)\ot E$.}\eea

 $V$ is a right comodule over $\T$, with the coaction $\delta_V=\db_V\ot\id_V:V\lora V\ot V^*\ot V$, hence a right comodule over $\E$. Since $\E$ is a bialgebra, $\Vn$ is also a right comodule over $\E$. Its dual, $(\Vn)^*$ is a left comodule over $\E$ in a canonical way, hence a right comodule over $\E^{\rm cop}$. The relation in \rref{eq3.1} implies that $R$ is a morphism of $\E$-comodules. Hence, $\Lambda_n$, $\S_n$  are right $\E$-comodules, they are factor comodules of $V^{\ot n}$. 
Actually, $\E_n$ is a subcoalgebra of $\E$ and the cocation of $\E$ on $\Vn,\Lambda_n,\S_n$ factorizes through $\E_n$. Since $\E_n$ is finite dimensional, ${\E_n}^*$ is an algebra and its left modules are in one to one correspondence with right $\E_n$-comodules. Therefore ${\E_n}^*$ acts on $\Vn,\Lambda_n,\S_n$. On the other hand, the operator $R$ induces a right action $\rho=\rho_R$ on $\Vn$ of the Hecke algebra $\H_n$. We have a quantum analogue of Schur's double centralizers theorem \cite[Theorem 2.1]{ph97}:
\bba \label{eq4} {\E_n}^*&\cong &\End_{\H_n}(\Vn),\\
\label{eq5}\rho(\H_n)&\cong&\End_{{\E_n}^*}(\Vn)=\End^{\E_n}(\Vn).\eea
As a consequence, simple $\E_n$ comodules are parameterized by primitive idempotents  of $\H_n$, which, in their order, are parameterized by partitions of $n$. For each primitive idempotents $E_\lambda$ of $\H_n$, $\rho(E_\lambda)$ is either zero of a simple $\E_n$-comodule, conjugated idempotents define isomorphic comodules. In particular, as $\E_n$-comodules
\bbs \Lambda_n\cong\Im\rho(Y_n),\quad \S_n\cong\Im\rho(X_n),\ees
where $X_n=([n]_q!)^{-1}\sum_{w\in\Ss_n}T_w$ , $Y_n=([n]_{q^{-1}}!)^{-1}\sum_{w\in\Ss_n}q^{-l(w)}T_w$, the $q$-symmetrizer and $q$-anti-symmetrizer operators, where $[n]_q:=(q^n-1)/(q-1)$. Moreover we have isomorphisms of algebras 
\bba\label{eq6} \S\cong\bigoplus_{n=0}^\infty\rho(X_n),\quad \Lambda\cong\bigoplus_{n=0}^\infty\rho(X_n),\eea
where the product on the spaces on the right-hand sides of these isomorphisms are given by
\bba\label{eq7} x\in\rho(X_n),y\in\rho(X_m)\lora x\cdot y:=\rho_{m+n}(X_{m+n})(x\ot y)\\
\label{eq8}
x\in\rho(Y_n),y\in\rho(Y_m)\lora x\cdot y:=\rho_{m+n}(Y_{m+n})(x\ot y).\eea
 This is the second realization of $\S$ and $\Lambda$, first considered by Gurevich \cite{gur1}. In other words, let $X=\sum_{n=0}^\infty\rho(X_n)$. Then $X$ is a projection on $\T(V)$, which carries the algebra structure of $\T(V)$ to its image $\Im X$ making this space an algebra. The second realization states that $\Im X$ is isomorphic to $\S_R$. Further, we have
\bba\label{eq9} \Im\rho(X_n)=\bigcap_{i=1}^{n-1}R(\Lambda)^n_i,\quad \Im\rho(Y_n)=\bigcap_{i=1}^{n-1}R(S)^n_i.\eea
This means $\S$ is isomorphic to $\Lambda^{!*}$, $\Lambda$ is isomorphic to $\S^{!*}$ as graded space, where $^*$ means graded dual.

\vspace{2ex}

\noindent{\it Example.} Let $R_d$ be Drinfel'd-Jimbo's solution of the Yang-Baxter equation of type $A_{d-1}$. Explicitly, $R_d$ is given as follows, with respect to a basis $x_1,x_2,\ldots,x_d$,
\bbas &&{R_d}^{kl}_{ij}=
\frac{q^2-q^{2\va_{ij}}}{1+q^{2\va_{ij}}}\delta_{ij}^{kl}+\frac{q^{\va_{ij}}(q^2+1)}{1+q^{2\va_{ij}}}\delta_{ij}^{lk},\quad 1\leq i,j,k,l\leq d,\quad \va_{ij}:=\mbox{sign }(j-i).\eeas
The Hecke equation of $R$ is $(R+1)(R-q^2)=0$. $\Im(R-q^2)$ has a basis consisting of elements of the form $x_i\ot x_j-q^{-1}x_j\ot x_i,i<j$.
$\Im(R+1)$ has a basis consisting of elements of the form $x_ix_j+qx_jx_i,i\leq j$. Thus the relations on the algebra $\S_{R_d}$ are $x_ix_j-q^{-1}x_jx_i=0,i<j$ and the relations on the algebra $\Lambda_{R_d}$ are $x_ix_j+qx_jx_i=0,i\leq j$. The second realization of $\S$ and $\Lambda$ says that their elements can be realized as $q$-symmetrized tensors in $k\langle x_1,x_2,\ldots,x_d\rangle $:
\bbas &&\mbox{for $\S$ }: x_{i_1} x_{i_2}\cdots x_{i_k}\loma \sum_{s\in \Ss_k}q^{-l(w)}x_{i_1s}\ot x_{i_2s}\ot\cdots\ot x_{i_ks}\\
&&\mbox{for $\Lambda$ }: x_{i_1}x_{i_2}\cdots x_{i_k}\loma \sum_{s\in \Ss_k}(-q)^{l(w)}x_{i_1s}\ot x_{i_2s}\ot\cdots\ot x_{i_ks}.\eeas

The situation for the algebra $\E$ is more complicated. In studying the Poincar\'e series of $\E$, Sudbery \cite{sud2} introduced an operator $\Phi$,  acting on $\E_1^{\ot 3}$, which is an analogue of the symmetrizer $X_n$. The study was further developed in \cite{ph97}, where a complete description of the Poincar\'e series of $\E$ was given in terms of the Poincar\'e series of $\Lambda$. Define for each $n\geq 0$ an operator
\bba\label{eq10}\Phi^n:=\sum_{w\in\Ss_n}(-1)^{l(w)}\bar R_w.\eea
As it has been mentioned, we identify $\bar R_w$ with ${\Rt^{-1}}_w\ot R_w$, acting on $\Vsn\ot \Vn$. For $n=3$, $\Phi^3$ is the operator introduced by Sudbery. It was shown that
\bba\label{eq12} \Im\Phi^n=\bigcap_{i=1}^{n-1}R(\E)^n_i,\eea
for transcendent $q$ \cite[Thm.~2.6]{ph97}, from which follows
\bba\label{eq13} P_\E(t)=P_\S\circ P_\Lambda,\eea
``$\circ$'' is the product in the $\lambda$-ring ${\Bbb C}_0[[t]]$ of power series with constant coefficient equal to 1 (see, e.g. \cite{knutson}).

\section{Another Realization for The Algebra $\E$}\label{sec3}
Unlike $X_n$ and $Y_n$, the operator $\Phi^n$ introduced in \rref{eq12} is not a projection. This is the main difficulty in giving the second realization for the algebra $\E$. If we know the eigenvalues of this operator, we can modify it to get a projection.

\subsection{The Operator $\Phi^n$}\label{sect-phi}
 Let us denote $R':=-qR^{-1}$. Then $R'$ also satisfies the equation $(x+1)(x-q)=0$; hence induces a representation of $\H_n$ on $\Vn$.  Let $\si=\si_R $ be the representation of $\H_n^{\rm op}$ on $(V^{\ot n})^*$: $\si (T_w)=(R'_w)^*$. If we identify $(V^{\ot n})^*$ with $V^{*\ot n}$ then $(R'_w)^* =\Rt'_{w^{-1}}.$ Applying $\si \ot \rho$ on the identity \rref{eq0}, we get
\bba\label{eq14}
\sum_{w\in\Ss_n}q^{-l(w)} \Rt'_{w}\ot R_w=\sum_{w\in\Ss_n}\si (E^{ij}_\lambda)\ot\rho(E^{ji}_\lambda).\eea
Recall that $\bar R=-q^{-1}{}^t{R'}\ot R$. Therefore,
\bba\nonumber
\sum_{w\in\Ss_n}(-1)^{l(w)}\bar R_w&=&\sum_{w\in\Ss_n}q^{l(w)}\Rt'_{w^{-1}}\ot R_w\\
&=&  \sum_{w\in\Ss_n}k^{-1}_\lambda\si (E^{ij}_\lam)\ot\rho(E^{ji}_\lam).\label{eq15}\eea

Denote 
$\displaystyle\Phi_\lam:=(\si \ot \rho)\Pi_\lam=\sum_{1\leq i,j\leq d_\lam}k^{-1}_\lam\si (E^{ij}_\lam)\ot\rho(E^{ji}_\lam).$
Then we have
\bba\label{eq16} \Phi_\lam^2=d_\lam k^{-1}_\lam\Phi_\lam.\eea
Thus, the spectrum of $\Phi^n$ is contained in $\{d_\lambda k^{-1}_\lam| \lam\part n\}\cup 0$. Note that if $q=1$, $d_\lam k^{-1}_\lam=n!.$

The isomorphisms in \rref{eq4} and \rref{eq5} imply the following decomposition of $\Vn$ as ${\E_n}^*-\H_n$-bimodule:
$$ \Vn\cong\bigoplus_{\lambda\part n} \Im\rho(E_\lambda)\ot S_\lam,$$
where $S_\lam$ is the simple right $\H_n$-module, isomorphic to the right ideal in $\H_n$, spanned by $\{ E_\lam^{jm}|1\leq m\leq d_\lam\}$, for each fixed $j$. Analogously, as ${\E_n}^{*\rm op}-\H_n^{\rm op}$-bimodules,
$$\Vns\cong\bigoplus_\lambda \Im\si (E_\lambda)\ot {S_\lam}^*,$$
where ${S_\lam}^*$ is the simple right $\H_n^{\rm op}$-module, isomorphic to the right ideal in $\H_n^{\rm op}$, spanned by $\{ E_\lam^{mj}|1\leq m\leq d_\lam\}$, for each fixed $j$. Therefore, 
\bba\label{eq17}\Vns\ot \Vn\cong \bigoplus_{\lambda,\mu\part n} \Im\si (E_\lambda)\ot \Im\rho(E_\mu)\ot {S_\lam}^*\ot S_\mu,\eea
as ${\E_n}^{*\rm op}\ot {\E_n}^*- \H_n^{\rm op}\ot\H_n$-bimodules. In particular, for any elements $x\in\Im\si(E_\mu), y\in\rho(E_\nu)$, $x\ot y\ot {S_\mu}^*\ot S_\nu$ is an invariant space of $\Phi^n$ and the action of $\Phi^n$ on this space does not depend on $x,y$. Moreover, we have:
\begin{lem}\label{lem11}Let $x\in \Im\si(E_\mu)$ and $y\in\Im\rho(E_\nu)$. Then the operator $\Phi_\lam$ vanishes on $x\ot y\ot {S_\mu}^*\ot S_\nu$, if $\mu\neq \lam$ or $\nu\neq \lam$; and if $\mu=\nu=\lam$, it has rank 1 on $x\ot y\ot {S_\lam}^*\ot S_\lam$.\end{lem}
\proof Since $S_\nu$ (resp. ${S_\mu}^*$) is simple over $\H_n$ (resp. $\H_n^{\rm op}$), the action of $\Phi^\lam$ on $x\ot y\ot {S_\mu}^*\ot S_\nu$ is equivalent to the action of $\Pi_\lam$. Thus, it is zero if $\mu\neq \lam$ or $\nu\neq \lam$. 

 Fix $i,j$ and let $\{ E_\lam^{mj}\ot E_\lam^{in} |1\leq m,n\leq d_\lam\}$ be a basis of ${S_\lam}^*\ot S_\lam$. We have, as elements in $\H_n^{\rm op}\ot \H_n$,
\bbas (E_\lam^{jm}\ot E_\lam^{in})\cdot \Pi_\lam&=&(E_\lam^{jm}\ot E_\lam^{in})\cdot \sum_{1\leq i,j\leq d_\lam}k^{-1}_\lam E^{ij}_\lam\ot E^{ji}_\lam\\
&=&\delta^m_n k^{-1}_\lam\sum_{1\leq l\leq d_\lam} E_\lam^{lj}\ot E_\lam^{il}\eeas
That is, the action of $\Pi_\lam$ from the right on ${S_\lam}^*\ot S_\lam$ has rank 1. 
\eee

\begin{cor}\label{cor1.1.1} As $\E^{\rm cop}\ot \E$-comodules
\bba\nonumber
\Im(\Phi_\lam)\cong \Im\si (E_\lam)\ot \Im\rho(E_\lam),\\
\Im(\Phi^n)\cong \bigoplus_{\lam\part n}\Im\si (E_\lam)\ot\Im\rho(E_\lam).\label{eq18}
\eea\end{cor}

Not all $\lam\part n$ contribute in the second decomposition of \rref{eq18}. For some $\lam$, $\si(E_\lam)$ or $\rho(E_\lam)$ my vanish. See Theorem \ref{lem-vanish} for the vanishing condition of $\si(E_\lam)$ and $\rho(E_\lam)$.

\subsection{The Operator $\Psi^n$.}\label{sect-psi}
The operator $\bar\Phi^n$ is not a right projection that we need for describing $\E$. The discussion in the previous subsection suggests us a new operator. According to \rref{eq3}, $\E_2$ is isomorphic to $\Im(\tilde R+q)$, where $\tilde R:=s_{(23)}\left(\Rt\ot R\right)s_{(23)}$. We shall also identify $\tilde R$ with $\Rt\ot R$ acting on $V^{*\ot n}\ot \Vn$. Let us define
\bbs \Psi^n:=\sum_{w\in\Ss_n}q^{-l(w)}\tilde R_w.\ees
 Let $\tau$ be the representation of $\H_n^{\rm op}$ on $(\Vn)^*$, induced by $R$. Thus, $\tau(T_w)=\Rt_{w^{-1}}$. In analogy to \rref{eq15}, we have
\bba\nonumber
\sum_{w\in\Ss_n}q^{l(w)}\tilde R_w= \sum_{w\in\Ss_n}k^{-1}_\lambda\tau (E^{ij}_\lam)\ot\rho(E^{ji}_\lam).\label{eq19}\eea
Notice that $\rho(E_\lambda)^*$, considered as $\H_n^{\rm op}$-module, is isomorphic to $\tau (E_\lambda)$. This is because $E_\lambda$ is also a primitive idempotent in $\H_n^{\rm op}$. Therefore, we have a decomposition of ${\E_n}^{*\rm op}\ot {\E_n}^*-\H_n^{\rm op}\ot\H_n$-bimodules:
\bba\label{eq191}\Vns\ot \Vn\cong\bigoplus_{\lam,\mu\part n}\Im\rho(E_\lam)^*\ot\Im\rho(E_\lam)\ot {S_\lam}^*\ot S_\lam.\eea
Set $\Psi_\lam=(\rho^*\ot \rho)\Pi_\lam.$ An analogue of Lemma \ref{lem11} holds for $\Psi_\lam$. Consequently, we have
\begin{cor}\label{cor12} There exist an isomorphism of $\E^{\rm   cop}\ot \E$-comodules
\bba\label{eq20}
\Im\Psi^n\cong\bigoplus_{\lam\part n}\Im\tau (E_\lam)\ot \Im\rho(E_\lam)\cong\bigoplus_{\lam\part n}\Im\rho(E_\lam)^*\ot \Im\rho(E_\lam).\eea\end{cor}

On the other hand, according to \rref{eq5}, there is an algebra isomorphism
\bbs {\E_n}^*\cong\bigoplus_{\lam\part n}\End_{\AA_\lam}(\Vn)\cong \bigoplus_{\lam\part n}\Im\rho(E_\lam)\ot \Im\rho(E_\lam)^*.\ees
Consequently, we have an isomorphism of coalgebras
\bba\label{eq21}\E_n\cong\bigoplus_{\lam\part n}\Im\rho(E_\lam)^*\ot \Im\rho(E_\lam).\eea
Therefore, as $\E^{\rm cop}_n\ot \E_n$-comodules,
\bba\label{eq22}\E_n\cong \Im\Psi^n.\eea
The operator $\Psi^n$ is not a projector. However we can slightly modify it to have a projection
\bbs \bar\Psi^n:=\sum_\lam k_\lam d^{-1}_\lam\tau (E_\lam^{ij})^*\ot\rho(E_\lam^{ji}).\ees
 Set $\Psi=\sum_{n=0}^\infty \bar\Psi^n.$
Then $\Psi$ is a projector on $\T(V^*\ot V)$, which then induces an algebra structure on its image $\Im\Psi$: for $a\in\Vsn\ot\Vn$ and $b\in V^{*\ot m}\ot V^{\ot m}$, $a\cdot b:=\bar\Psi^{n+m}(a\ot b)$.
\begin{thm}\label{thm1} Assume that the parameter $q$ is not a root of unity of order greater that 1. Then the projection $\Psi$ induces an algebra isomorphism from  $\E=$ to $\Im\Psi$.\end{thm}
\proof We have $(\tilde R+q)(\bar R-1)=(\bar R-1)(\tilde R+q)=0$. Therefore, if $x\in\Im(\bar R_i-1)$ then $(\tilde R+q)x=0$, consequently, $\Psi^n$ vanishes on $\sum_{i=1}^{n-1}R(\E)_i^n$; hence $\bar \Psi^n$ vanishes on $\sum_{i=1}^{n-1}R(\E)_i^n$, too. Taking \rref{eq22} into account, we conclude that $\sum_{i=1}^{n-1}R(\E)_i^n$ is precisely the Kernel of $\bar\Psi^n$. Thus, the linear map $\Psi:\T(V^*\ot V)\lora\bigoplus_{n=0}^\infty\Im\bar\Psi^n$, $a\loma \oplus_{n=0}^\infty\bar\Psi^n(a)$ has the following properties
\bba\bbar{l}a\in R(\E)\Longrightarrow \Psi(a)=0\\
a=\bar\Psi(a) \bmod  R(\E).\eear\label{eq23}
\eea
That is $R(\E)=\Ker\Psi$.
Therefore, $\Psi$ factorizes through $R(\E)$ to a linear isomorphism $\psi:\E\lora \bigoplus_{n=0}^\infty\Im\bar\Psi^n$.

Let $a\in(V^*\ot V)^{\ot m}$, $b\in(V^*\ot V)^{\ot n}$. According to \rref{eq23}, we have
\bbs \bar\Psi^m(a)=a\bmod R(\E)^m\quad  \bar\Psi^n(b)=a\bmod R(\E)^n,\ees
hence
\bbs \bar\Psi^m(a)\ot \bar\Psi^n(b)=a\ot b\bmod R(\E)^{m+n}.\ees
Consequently,
\bba
\bar\Psi^{m+n}(\bar\Psi^m(a)\ot\bar\Psi^n(b))=\bar\Psi^{m+n}(a\ot b)\label{eq24},\eea
meaning that $\psi$ is an algebra homomorphism and therefore isomorphism $\E\lora \Im\Psi$.\eee

We now proceed to show that, for $q$ not a root of unity of order greater than 1 (cf. Eq. \rref{eq12}),
\bba\label{eq25} \Im\Phi^n=\bigcap_{i=1}^{n=1}R(\E)^n_i.\eea
The inclusion ``$\subset$'' is obviously, for we have $\Phi^n=(\bar R_i-1)P_i$ for certain operator $P_i$, $i=1,2,\ldots,n-1$. 
To show the equality, we compare the dimensions of ${\E^!}_n$ and $\Im\Phi^n$. Let $l_\lambda=\dim_k\Im\rho(E_\lambda)$. From the definition of the operator $S$, we see that $\si (E_\lam)\cong \rho(E_{\lam'})$ as vector spaces. Therefore $\dim\Im\Phi^n=\sum_{\lam\part n}l_{\lam'}l_\lam$.

 On the other hand, since $\E$ is a Koszul algebra (with the above assumption on $q$), (cf. \cite[Thm~2.5]{ph97})
\bba\label{eq26} P_\E(t)P_{\E^!}(-t)=1.\eea
Hence, according to \rref{eq13}, we have $P_{\E^!}(t)=P_\Lambda(t)\circ P_\Lambda(t)$, that is dim${\E^!}_n=\sum_{\lam\part n} l^2_\lam,$ \cite{ph97c}. Therefore 
\bba \dim\bigcap_{i=1}^{n-1}R(\E)^n_i=\dim ({{\E^!}_n}^*)=\dim E^!_n=\sum_{\lam\part n}l_{\lam'}l_\lam.\eea
Thus, $\dim\Im\Phi^n=\dim {\E^!}_n$. Whence \rref{eq25} follows.

Finally, let us denote by $\F$ the quadratic algebra on $V^*\ot V$ with relation $R(\F):=\Im(\tilde R+q)$. Thus, $\F$ can be considered as the quantum exterior algebra over the matrix quantum semi-group. We have
\bba\label{eq28} R(\F)=s_{(23)}\left(R(\S)^\perp\ot R(\Lambda)\bigoplus R(\Lambda)^\perp\ot R(\S)\right).\eea
In terms of the matrix $E$, the relations can be given as follows:
\bbs RE_1E_2=-qE_1E_2 R.\ees

The vector space $V^*\ot V$ is self dual with respect to the pairing $(\phi\ot x,\psi\ot y):=(\phi,y)(\psi,x)$. With respect to this pairing, $\F$ is canonically isomorphic to $\E^!$. Therefore $\F$ is Koszul algebra and
\bba\label{eq29}P_\F(t)=P_\E(-t)^{-1}=P_\Lambda(t)\circ P_\S(t),\eea
or equivalently
\bba\label{eq30} \dim \F_n=\sum_{\lam\part n}l_\lam l_{\lam'}.\eea

Decompose $(V^*\ot V)^{\ot n}$ (which is identified with $\Vsn\ot \Vn$) into simple $\H_n^{\rm op}\ot \H_n$-modules (with the action given by $\tau \ot\rho$). Notice that, as in the case of $\E$, $\F_n$ and $R(\F)^n$ are $\H_n^{\rm op}\ot \H_n$-modules. Recall that the operator $\Phi^n$ also has rank $\sum_{\lam\part n}l_\lam l_{\lam'}$ on $\Vsn\ot \Vn$ and vanishes on $R(\F)^n$. Comparing the dimension we see that $R(\F)^n$ is precisely the kernel of $\Phi^n$. Analogously, the image of $\Psi^n$ is $\bigcap_{i=1}^{n-1}R(\F)^n_i$. Thus we proved an analogue of Theorem \ref{thm1}:
\begin{thm}\label{thm2} The algebra $\F$ is a Koszul algebra and can be realized as $\Im\Phi$, $\Phi:=\bigoplus_{n=0}^\infty\bar\Phi^n$.\end{thm}

The question when $\rho(E_\lam)$ is zero can be answered by knowing the Poincar\'e series of $\S_R$. More precisely, is is proved that the Poincar\'e series of $\S_R$ is a rational function having only negative roots and positive pole (as a complex function). Let $r$ denote the number of poles and $s$ denote the number of roots of $P_\S(t)$. We call the pair $(r,s)$ the {\it ~birank} of $R$. For example, the birank of the operator $R_d$ in Section \ref{sec1} is $(d,0)$. 

Let $\Gamma_{r,s}$ denote the set of partitions $\lam$ such that $\lam_{r+1}\leq s$.
\begin{thm}\cite[Theorem 5.1]{ph97c}\label{lem-vanish}
Assume that the operator $R$ has the birank $(r,s)$. Then the comodule $V_\lam:=\Im\rho(E_\lam)$ is non-zero if and only if $\lam\in\Gamma_{r,s}$.\end{thm}

\section{The Quantum Spaces of Homomorphisms}\label{sec4}
The notion of $\E_R$ as an ``endomorphism ring'' of a quantum space can be generalized to the notion of ``space of homomorphisms'' of two quantum spaces.
Let $R$ and $S$ be Hecke operators on $V$ and $W$, respectively. We define the quadratic algebra $\M=\M_{SR}$ on $W^*\ot V$, whose relation is
\bba\label{eq31}R(\M):=s_{(23)}\left((R(\S_S)^\perp\ot R(\S_R)\oplus R(\Lambda_S)^\perp\ot R(\Lambda_R)\right),\eea
where, as usual, $R(\S_S):=\Im(S-q)$, $R(\S_R):=\Im(R-q)$, and so on. We have
\bba\label{eq32} R(\M)=\Im(s_{(23)}(\St^{-1}\ot R)s_{(23)}-\id).\eea
As in the previous section, we shall identify the two vector space $(W^*\ot V)^{\ot n}$ and $W^{*\ot n}\ot \Vn$.

The algebra $\M$ can be interpreted as the function algebra on the quantum space of homomorphisms (or quantum hom-space) from the quantum space associated to $R$ to the one associated to $S$. Let $x_1,\cdots,x_m$ be a basis of $V$ and $\eta^1,\cdots,\eta^n$ be a basis of $W^*$. Then $\{m^i_j:=\eta^i\ot x_j\}$ form a basis of $W^*\ot V$. $\M$ is then isomorphic to
\bba\label{eq33}
\bK\langle m^i_j|1\leq i\leq n,1\leq j\leq m\rangle/(SM_1M_2-M_1M_2 R)\eea
where $M=(m^i_j)$, $M_1:=M\ot \id(n)$, $M_2:=\id(m)\ot M$. $\M$ has the following properties: Let $A:=(a_1,a_2,\ldots,a_n)$ be a point of $\S_S$, i.e., $(A\ot A)(R-q)=0$, and $T=(t^i_j)$ be a point of $\M_{SR}$, such that $a_i$ and $t_l^k$ commute. Then $A\dot\ot M$ is a point of $\S_R$, where $(A\dot\ot M)_i:=\sum_k a_k\ot m^k_i$. Analogously, if $B$ is a point of $\Lambda_S$ commuting with $T$ then $B\dot\ot T$ is a point of $\Lambda_R$.

There is also an interpretation of $\M_{SR}$ from the categorical view-point. The bialgebra $\E$ can be constructed as the Coend of the functor $\Ff$ from the braided monoidal category $\V$, generated by one object $v$ and one morphism $\tau:v^{\ot 2}\lora v^{\ot 2}$, into the category of vector space, such that $\Ff(v)=V$ and $\Ff(\tau)=R$ (cf. \cite{schauen,majid1}). That is, for any vector space $X$,
\bba\label{eq34}\Nat(\Ff,\Ff\ot X)\cong\Hom_\bK(\E,X)\eea
where $\Nat(\Ff,\Gg)$ denotes the set of natural transformation between functors $\Ff$ and $\Gg$, $\Ff\ot X$ is the functor that sends $w$ to $\Ff(w)\ot X$ and sends $f$ to $\Ff(f)\ot\id_X$, $v\in \V, f\in\Mor(\V)$.
Let us now consider another functor $\Gg$, with $\Gg(v)=W$ and $\Gg(\tau)=T$. Then we have
\bba\label{eq35}\Nat(\Ff,\Gg\ot X)\cong \Hom(\M_{SR},X).\eea

The exterior algebra on the quantum hom-space is defined to be
\bba \N=\N_{SR}:=\T(W^*\ot V)/\Im(s_{(23)}(\St\ot R)s_{(23)}+q\cdot \id).\eea

The bialgebra $\E_R$ coacts on $\M$ and $\N$ from the right. The coaction is induced from the one on $W^*\ot V$: $\delta(m^i_j)=\sum_km^i_k\ot e^k_j.$ Analogously, $\E_S$ coacts on $\M$ and $\N$ from the left, with the coaction induced from $\delta(m^i_j)=\sum_le^i_l\ot m^l_j.$ Thus, $\M$ and $\N$ are right $\E_S^{\rm cop}\ot \E_R$-comodule algebras.

We show in this section that $\M$ and $\N$ are Koszul algebras, compute their Poincar\'e series and give a second realization. Since $N_{SR}=\M_{S'R}$, $S'=-qS^{-1}$, it is sufficient to study $\M_{SR}$.

As in \rref{eq17}, we have a decomposition of $(W^*\ot V)^{\ot n}$ as an ${\E_S}^{\rm cop}\ot \E_R- \H_n^{\rm op}\ot\H_n$-bimodule
\bba\label{eq36}W^{*\ot n}\ot \Vn\cong \bigoplus_{\lambda,\mu\part n} \Im\si(E_\lambda)\ot \Im\rho(E_\mu)\ot {S_\lam}^*\ot S_\mu.\eea
The subspace $R(\M)^n$ of $(W^*\ot V)^{\ot n}$ is an ${\E_S}^{\rm cop}\ot \E_R$-comodule, hence
\bbs R(\M)^n=\bigoplus_{\lambda,\mu\part n}\left(R(\M)^n\cap \Im\si_S(E_\lambda)\ot \Im\rho_R(E_\mu)\ot {S_\lam}^*\ot S_\mu\right).\ees

\vskip1ex

\noindent{\it Remark.}
if the action of $\St_i\ot R_i$ on ${S_\lam}^*\ot S_\mu$ is not zero then this action does not  depend on $S$ and $R$. In fact, the action of $\St_i\ot R_i$ on ${S_\lam}^*\ot S_\mu$ in this case is the action of $T_i\ot T_i$. 

 Define the operators $\Phi_{SR}$ and $\Psi_{SR}$ as in \rref{eq10} with $\bar R=q\St^{-1}\ot R$ and $\bar R=\St\ot R$, respectively. The corresponding projectors $\bar\Phi_{SR}$ and $\bar \Psi_{SR}$ are defined similarly. Notice that $\Phi_{SR}=\Psi_{S'R}$. From the proof of Theorem \ref{thm1} and using the above remark we have
$$ R(\M)^n=\Ker\bar\Psi_{SR}^n.$$
Consider the action of $\bar\Psi_{SR}^n$ on a module $x\ot y\ot {S_\lam}^*\ot S_\mu$, $x\in\si(E_\lam), y\in\rho(E_\mu)$. This has rank 1 if $\lam=\mu$ and 0 otherwise. Therefore, for $k_\lambda:=\dim\Im\rho_S(E_\lambda), l_\lambda:=\dim\Im\rho_R(E_\lambda)$, we have
\bbs \dim_k\M_n=\dim_k(\M_1^{\ot n}/ R(\M)^n)=\sum_{\lambda\part n}l_\lambda k_\lambda.\ees

Since $\bar\Psi_{SR}^n$ is a projector, so is the map
$\Psi_{SR}:\bigoplus_{n=0}^\infty\bar\Psi^n_{SR}$, which induces an isomorphism of algebras $\M\lora \bigoplus_{n=0}^\infty\Im\Psi_{SR}$. This map is also a homomorphism of $\E_S^{\rm cop}\ot \E_R$-comodules because each $\bar\Psi_{SR}$ is.

 The remark above implies that the lattice induced by $R(\M)^n_i\cap({S_\lam}^*\ot S_\mu),1\leq i\leq n-1$ is distributive in $({S_\lam}^*\ot S_\mu)$ (see \cite{ph97} for more details). Consequently, the lattice generated by $R(\M)^n_i,1\leq i\leq n-1$ is distributive in $W^{*\ot n}\ot \Vn$, that is $\M_{SR}$ is a Koszul algebra.

\begin{thm}Let $R=R_q$ and $S=S_q$ be Hecke operators, where $q\in\bK^\times$ is not a root of unity of order greater than 1. Then the algebra $\M_{SR}$ is a Koszul algebra, its Poincar\'e series are given by
\bba P_{\M}(t)=P_{\S_R}\circ P_{\S_S}(t).\eea
There is a realization of $\M$ a subspace of $\T(W^*\ot V)$: the following is an isomorphism of $\E_S-\E_R$-bicomodule algebras
\bba \M\cong \bigoplus_{n=0}^\infty\Im\bar\Phi_{SR}^n\eea
where, for each $n$, as $\E_S-\E_R$-bicomodules
\bbs\M_n\cong\Im\bar\Psi_{SR}^n=\bigoplus_{\lam\part n}\Im\rho_S(E_\lam)^*\ot\Im\rho_R(E_\lam).\ees
Further, we have
\bba 
\Im(\bar\Psi^n_{SR})=\bigcap_{i=1}^{n-1}R(\M)^n_i.\eea
\end{thm}

\section{Quantum Determinantal Ideals}\label{sec5}

In commutative algebra, the ideal $I_k$, generated by the $k\times k$-minors in the coordinate ring of the varieties $M_{k}{(m,n)}$ is called determinantal ideal ($0\leq k\leq\min(m,n)$). This ideal is invariant with respect to a natural action of the group $G=GL_k(m)\times GL_k(n)$ on $\O(M_k(m,n))$. It is proved to be prime, see, e.g. \cite{cep1}. The variety determined by $I_k$ is called determinantal variety.

In the quantum setting, we have a coaction of the bialgebra $\G=\E^{\rm cop}_S\otimes \E_R$ on the algebra $\M_{SR}$. A subspace (resp. two-sided ideal) in $\M_{SR}$, which is invariant with respect to the coaction of $\G$ will be called invariant subspace (resp. invariant ideal). Let $M_\lambda$ denote $\Im\tau(E_\lambda)\ot \Im\rho(E_\lambda)$. Then any invariant subspace of $\M$ is a direct sum of some $M_\lambda$, $\lam\in\P$, ($\P$ is the set of all partitions).

Let us denote by $\P_{SR}$ the set of partitions $\lambda$, such that $k_\lambda l_\lambda\neq 0$ (see the previous section). This set can be fully described using Theorem \ref{lem-vanish}. In this section we show that there is a one-one correspondence between $\G$-invariants ideals in $\M$ and D-ideals in $\P_{SR}$. The latter is defined as follows: a subset $J$ of $\P_{SR}$ is called a D-ideal if for any $\si\in J$ and any $\tau\in\P_{SR}$, such that $\tau\supset \si$, one has $\tau\in J$ (cf. \cite{cep1}). For any subset $J$ of $\P_{SR}$, let $I(J)$ denote the subspace $\sum_{\si\in J}M_\si$.

We need the following key lemma. 
\begin{lem}\label{key-lem} Let  $\lambda$ and $\mu$ be partitions, such that all the Littlewood-Richardson coefficients $c_{\lam\mu}^\gamma$ is at most 1. Let $C_{\lambda,\mu}$ denote the set of partitions $\gamma$ such that $c_{\lambda,\mu}^\gamma=1$. Then the image of the product of $\M$, restricted on $M_\lambda\ot M_\mu$, is $I(C_{\lambda,\mu})$.\end{lem}
Reference for the Littlewood-Richardson coefficients is \cite{mcd2}. The proof of this lemma will be given at the end of this section.

Notice that for any partition $\lam$, $\lam$ and $(1^k)$ satisfies the condition of the lemma above. In fact, using the Littlewood-Richardson rule $c_{\lam(1^k)}^\gamma$ is equal to 1 if and only if $\gamma\supset\lam$ and $\gamma_j-\lam_i\leq 1$, $|\gamma|-|\lam|=k$; otherwise $c_{\lam(1^k)}^\gamma=0.$  Thus we have
\begin{cor}\label{cor-mult} For any partition $\lam$ and any integer $k$, $M_\lam\cdot M_{(1^k)}=I(C_{\lam,(1^k)})$.\end{cor}
Following \cite{cep1}, we denote by $I_\si$ the ideal in $\M$, generated by $M_\si$. We have
\begin{thm}\label{thm3} $I_\si=\sum_{\tau\supset \si}M_\tau.$\end{thm}
\proof
Since $M_\si\cdot M_{(1)}=I(C_{\si(1)})$ and $\gamma\in C_{\si(1)}$ implies $\gamma\supset \si$, we have $I_\si\subset \sum_{\tau\supset \si}\M_\tau.$

Assume $\tau\supset \si$ and $|\tau|=|\si|+1$. Then $\tau\in C_{\lam(1)}$. Hence, $M_\tau$ appears in $M_\si\cdot M_{(1)}$, thus $M_\tau\subset I_\si$. By induction, $\sum_{\tau\supset \si}M_\tau\subset I_\si.$\eee

For each set $J\subset \P$, let $\langle J\rangle $ denote the smallest D-ideal in $\P$, containing $J$. The following facts follow immediately from Theorem \ref{thm3}.

\begin{cor}\begin{enumerate}\item The ideal generated by $I(S)$ for some subset $S\subset \P$ is $I(\langle S\rangle )$.
\item $I_\si\supset I_\tau$ if and only if $\tau\supset \si$.
\item Let $J\subset\P$. Then $I(J)$ is an ideal in $\M$ if and only if  $J$ is a D-ideal in $\P$.\end{enumerate}\end{cor}

We now describe the product of determinantal ideals. We need the following order on the set of partitions (cf. \cite{cep1}). Firstly, we identify a partition with the diagram it determines. Define functions $\beta_k$, $k=1,2,...$, on the set of diagram as follows. For a partition $\si=(\si_1,\si_2,\cdots)$, $\beta_k(\si):=\si'_1+\si'_2+\cdots+\si'_k$, where $\si'$ is the conjugate partition of $\si$. In other words, $\beta_k(\si)$ is the number of boxes in the first $k$ columns of the diagram determined by $\si$. We say that $\tau\geq \si$ if for all $k$, $\beta_k(\tau)\geq \beta_k(\si)$.
\begin{pro}
Let $\si=(\si_1,\si_2,...)$ be a partition. Then the product of ideals $I_{\si_i}$ is $I(D_{\sigma'})$ where $D_{\si'}$ is the set of partitions $\tau$ such that $\tau\geq \si'$ in the above defined order.\end{pro}
\proof We use induction. If $\si$ has only one non-zero component, the claim is obvious. Let $\si=(\si_1,\si_2,\ldots,\si_r)$. Denote $\si_*$ the partition $\si=(\si_1,\si_2,\ldots,\si_{r-1})$. Using the induction assumption, one reduces it to showing that
$$I(D_\si)=I(D_{\si_*})\cdot I_{\si_r}.$$

Let $M_\lam\subset I( D_{\si_*})\cdot I_{\si_r}.$ Then $\lam\in C_{\gamma(1^{\si_r})}$, for some $\gamma\in D_{\si_*}$, by Corollary \ref{cor-mult}. In this case $\lam$ contains $\gamma$ and satisfies $\lam_i-\gamma_i\leq 1$, $|\lam|-|\gamma|=\si_r.$ Hence, for any $k=1,2,...,r-1$, $\beta_k(\lam)\geq \beta_k(\gamma)\geq\beta_k(\si_*)=\beta_k(\si)$, and for $k=r$, $\beta_r(\lam)=\beta_{r-1}(\gamma)+\si_r\geq \beta_{r-1}(\si_*)+\si_r=\beta_r(\si)$. Thus, $\lam\in\D_\si$, that is $I(D_{\si_*})\cdot I_{\si_r}\subset I(D_\si).$

To show the converse inclusion we need the following lemma.
\begin{lem}
 With the above notations, if $\lam\geq \si'$ then there exists $\gamma\subset \lam$, such that $\lam\geq {\si_*}'$ and $\lam\in C_{\gamma(1^{\si_r})}$.\end{lem}
\proof  Instead of considering partitions $\lam$ and $\si$, we consider their diagrams using the same notation. By assumption, the number of boxes in the first $k$ columns of $\lam$ is greater than the corresponding number in $\si$, for $k=1,2,...$ Let $b(\lam)$ denote the set of boxes lying in the rightest place in each row. Thus the number of boxes in $b(\lam)$ is $\beta_1(\lam)$. Since $\lam\geq \si'$, $|b(\lam)|\geq b_1(\si')=\si_1\geq \si_r$. To obtain $\gamma$, we remove from $\lam$ $\si_r$ boxes in $b(\lam)$. Starting from the lowest one in the rightest column we remove all the boxes bottom up (all these boxes belong to $b(\lam)$). Then we remove those boxes of $b(\lam)$ lying in the second rightest column in the same way and keep doing this further to the left. We remove as many boxes as the number of boxes in the rightest column of $\si'$, i.e. $\si_r$ boxes. The diagram obtained is $\gamma$. 

To see that $\gamma\geq {\si_*}'$, we proceed as follows. Each time when we remove a box from $\lam$, we also remove a box from the rightest column of $\si'$ in the order bottom up. Thus we have two sequences of diagram denoted by $\lam=\lam(0),\lam(1),\ldots,\lam(\si_r)=\gamma$ and $\si'=\si'(0),\si'(1),\ldots,\si'(\si_r)={\si_*}'$. It is easy to see that if $\lam(t)\geq \si'(t)$ then $\lam(t+1)\geq \si'(t+1)$. Since $\lam\geq \si'$, we conclude that $\gamma\geq{\si_*}'$. 

On the other hand the set $C_{\gamma(1^{\si_r})}$ contains such partitions $\lam$ that $0\leq \tau_k-\gamma_k\leq 1$ and $|\tau|-|\gamma|=\si_r$. Thus, the construction of $\gamma$ above also implies that $\lam\in C_{\gamma(1^{\si_r})}$. The Lemma is therefore proved.

Assume now that $\lam\in D_\si$, by Lemma, there exist  $\gamma\subset\lam$ such that $\gamma\geq\si_*$ and $\lam\in C_{\gamma(1^{\si_r})}$.  That is $\gamma\in D_{\si_*}$ and, by Corollary \ref{cor-mult}, $M_\lam\subset I(D_{{\si_*}'})\cdot I_{\si_r}$. That is,  $I(D_{\si'})\subset I(D_{\si_*})\cdot I_{\si_r}.$
The proposition is proved.\eee

The next interesting question is whether the primeness, radicals, preserve through the correspondence $J\lora I(J)$. This problem is completely open. The following property is the only one in this direction that I am able to prove.
\begin{pro} Let $J$ be an invariant ideal in $\M$. Then $\sqrt{J}$ is invariant, too.\end{pro}
We first recall the definition of the radical of a (two-sided) ideal in a non-commutative ring. An $m$-system in a ring $R$ is a set of elements such that for each two of its elements $a,b$, there exists an element $r$ from the ring, such that $arb$ belong to this set. Let $I$ be an ideal. The radical $\sqrt{I}$ of $I$ is defined to be the set of elements $s$, such that any $m$-system containing $s$ intersects non-trivially with $I$.

\proof We have to show that for any element $a\in\sqrt{I}$, the $\G$-invariant space generated by $a$ also belong to $\sqrt{I}$. To simplify the discussion, we introduce the algebra $\G^*$. Recall that $\G=\E_S^{\rm op}\ot \E_R$ is a cosemisimple bialgebra. Hence its dual is the completion of a direct sum of endomorphism rings on certain vector spaces
$$\G^*=\overline{\bigoplus_{\lam\mu}\End(V_{\lam\mu})}\cong\prod_{\lam\mu}\End(V_{\lam\mu}),$$
where the product on the right-hand side is defined componentwise. Let $\pi_{\lam\mu}$ be the projection on the $\lam\mu$-component. It is obvious that an element $f$ of $\G^*$ is invertible if and only if $\pi_{\lam\mu}f$ is invertible in $\End(V_{\lam\mu})$, for all $\lam,\mu$. Let $U$ be the of invertible elements in $\G^*$.
Since the set of invertible elements in $\End(V_{\lam\mu})$ spans this vector space, the set $U\cdot m$ spans the submodule generated by $m$. Let now $m\in\sqrt{I}$. It is thus sufficient to check that $U\cdot m\in\sqrt{I}$. This is an easy consequence of the fact that $U$ consists of invertible elements and the product on $\M$ is $\G^*$-equivariant. The proposition is proved.\eee

\noindent{\it Proof of Lemma \ref{key-lem}.} We use the characterization of the product of $\M$ given in Theorem \ref{thm1}. Thus for $a\in M_\lambda$, $b\in M_\mu$, $a\cdot b=\bar\Psi^{l+m}(a\ot b)$, where $l=|\lam|, m=|\mu|$. It is to show that for any given $\gamma\in C_{\lambda\mu}$ there exist such $a$ and $b$, that $\bar\Psi^\gamma(a\ot b)\neq 0$ . According to the fact that the product is $\G$-equivariant and that $M_\sigma$ is a simple $\G$-comodule the assertion of Lemma will then follow.

Denote $V_\lam:=\Im\rho(E_\lam)$. We have an isomorphism of ${\E_n}^{*\rm op}\ot {\E_n}^*\ot \H_n^{\rm op}\ot \H_n$ modules (cf. Eq. \rref{eq191})
$$V^{*\ot n}\ot \Vn\cong V_\lam^*\ot V_\lam\ot {S_\lam}^*\ot S_\mu.$$

For any fixed $i$ and $j$, $1\leq i,j\leq d_\lam$, we can assume that the set $\{ E_\lam^{mi}|1\leq m\leq d_\lam\}$ and $\{ E_\lam^{jm}|1\leq m\leq d_\lam\}$ are bases of ${S_\lam}^*$ and $S_\lam$, respectively. In this setting, an element of $M_\lambda$ can be represented as a linear combination of elements of the form $v\ot(\sum_{m}E_\lam^{mi}\ot E_\lam^{jm})$ for some $v\in V^*_\lam\ot V_\lam$. That is, as a subspace of $V_\lam^*\ot V_\lam\ot {S_\lam}^*\ot S_\lam$,
$$M_\lam=V_\lam^*\ot V_\lam\ot (\sum_{m}E_\lam^{mi}\ot E_\lam^{jm}).$$

We have, with the assumption on $\lam$ and $\mu$
 \bba\label{}V_\lam\ot V_\mu\cong \bigoplus_{\si\in C_{\lam,\mu}}V_\si,\quad V^*_\lam\ot V^*_\mu\cong \bigoplus_{\si\in C_{\lam,\mu}}V^*_\si,\eea
as $\E_{l+m}$-comodules, where the coaction of $\E_{l+m}$ on the left-hand side is induced by the product on $\E$, $\E_l\ot \E_m\lora \E_{m+l}$. Therefore,
$$M_\lam\ot M_\mu= \bigoplus_{\si,\eta\in C_{\lam,\mu}}V^*_\si\ot V_\eta\ot\left( \sum_{m,n}E_\lam^{mi}\ot E_\lam^{jm}\ot E_\mu^{nl}\ot E_\mu^{kn}\right),$$
so that, for each $\gamma\in C_{\lam,\mu}$, we can choose $v\in V^*_\lam\ot V_\lam$ and $w\in V^*_\mu\ot V_\mu$, such that the projection $\pi_\gamma(v\ot w)$ of $v\ot w$ on $V^*_\gamma\ot V_\gamma$ through the above isomorphism is not zero.  The crucial point here is that in the decomposition \rref{} is multiplicity-free.

On the other hand, if we embed $\H_l\ot \H_m$ into $\H_{l+m}$ in the standard way, then $S_\gamma$, considered as $\H_l\ot \H_m$-module, contains $S_\lam\ot S_\mu$ as a simple subcomodule. Therefore, the space
$$V^*_\gamma\ot V_\gamma\ot \left(\sum_{m,n}E_\lam^{mi}\ot E_\lam^{jm}\ot E_\mu^{nl}\ot E_\mu^{kn}\right)$$
is a subspace of $V^*_\gamma\ot V_\gamma\ot {S_\gamma}^*\ot S_\gamma$. The element $(a\ot b)\bar\Psi^\gamma$ is then
$$\pi_\gamma(v\ot w)\ot\left(\sum_{m,n}E_\lam^{mi}\ot E_\lam^{jm}\ot E_\mu^{nl}\ot E_\mu^{kn}\right)\bar\Psi_\gamma.$$
Here, the element $(\sum_{m}E_\lam^{mi}\ot E_\lam^{jm})$ is considered as an element of $\H_{l+m}$ by the embedding $\H_l\ot\H_m\lora \H_{l+m}$. Thus, it remains to show that (in $\H_n^\op\ot\H_n$)  
\bba\label{eq39}\left(\sum_{m,n}E_\lam^{mi}\ot E_\lam^{jm}\ot E_\mu^{nl}\ot E_\mu^{kn}\right)\Pi_\gamma\neq 0.\eea

Since the sets $\{E_\lam^{jm}, 1\leq m\leq d_\lam\}$, for $i=i_1$ and $i=i_2$ can be obtained from each other by multiplying with $E_\lam^{i_1i_2}$ on the left, \rref{eq39} does not depend on $i,j$ and $k,l$. That means, if \rref{eq39} holds (or fails) for some $i,j,k,l$, it should holds (or fails) for all $i,j,k,l$. So that, if on the left-hand side of \rref{eq39} we set $i=j$, $k=l$ and summing it up after these indices and show that this element is not zero, we will be done. Thus, we have to show that
\bbas\left(\sum_{mi}E_\lam^{m,i}\ot E_\lam^{im}\ot E_\mu^{nk}\ot E_\mu^{kn}\right)\Pi_\gamma\neq 0,\eeas
or $\Pi_\lam \Pi_\mu\Pi_\gamma\neq 0$,
for any $\gamma\in C_{\lam,\mu}$, here, $\Pi_\lam$ and $\Pi_\mu$ are considered as elements of $\H_{l+m}^{\rm op}\ot \H_{l+m}$. Let $F_\lam$ be the minimal central primitive element in $\H_l$, corresponding to $\lam$. Then $F_\lam$ is also a central primitive element in $\H_l^{\rm op}$. Hence from the definition of $\Pi_\lam$ and \rref{eq0}, we have
\bba\label{eq40} \Pi_\lam=k_\lam^{-1}(F_\lam\ot F_\lam)\left(\sum_{w\in\Ss_l}q^{-l(w)}R_{w}\ot R_{w^{-1}}\right).\eea
Analogous equalities hold for $\Pi_\mu$ and $\Pi_\gamma$. Let $\D_\lam$ be the set of left coset representatives of $\Ss_l$ in $\Ss_{l+m}$, such that $l(w)l(t)=l(wt)$ for $w\in\Ss_l,t\in\D_\lam$. Then we have, as elements of $\H_{l+m}^{\rm op}\ot \H_{l+m}$,
$$\sum_{v\in\Ss_{l+m}}q^{-l(v)}R_{v}\ot R_{v^{-1}}= \left(\sum_{w\in\Ss_l}q^{-l(w)}R_{w}\ot R_{w^{-1}}\right)\left(\sum_{t\in\D_\lam}q^{-l(t)}R_{t}\ot R_{t^{-1}}\right).$$
According to \rref{eq40}, we have
 \bbs\Pi_\lam\sum_{w\in\Ss_l}q^{-l(w)}R_{w}\ot R_{w^{-1}}=d_\lam k_\lam^{-1}(F_\lam\ot F_\lam)\sum_{w\in\Ss_l}q^{-l(w)}R_{w}\ot R_{w^{-1}}\ees
Therefore, 
\bbas\lefteqn{\Pi_\lam\Pi_\gamma=\Pi_\lam\left(\sum_{w\in\Ss_l}q^{-l(w)}R_{w}\ot R_{w^{-1}}\right)\left(\sum_{t\in\D_\lam}q^{-l(t)}R_{t}\ot R_{t^{-1}}\right)(F_\gamma\ot F_\gamma)}\\
&=&d_\lam k_\lam^{-1}(F_\lam\ot F_\lam)\left(\sum_{w\in\Ss_l}q^{-l(w)}R_{w}\ot R_{w^{-1}}\right)\left(\sum_{t\in\D_\lam}q^{-l(t)}R_{t}\ot R_{t^{-1}}\right)(F_\gamma\ot F_\gamma)\\
&=& d_\lam k_\lam^{-1}(F_\lam\ot F_\lam)(F_\gamma\ot F_\gamma)\sum_{w\in\Ss_{l+m}}q^{-l(w)}R_{w}\ot R_{w^{-1}}.\eeas

Thus, it is led to showing that $F_\lambda F_\mu F_\gamma\neq 0$. This is obvious by the assumption, that $\gamma\in C_{\lam,\mu}$. The Lemma is therefore proved.\eee

\section{Invariant theory}\label{sec6}
Let $m,n,t$ be positive integers. The group $GL(t)$ acts on the variety $M(m,t)\times M(t,n)$ in the following way:
\bba\label{eqinv1} g(A,B)=(Ag^{-1},gB), g\in GL(t), A\in M(m,t), B\in M(t,n).\eea
This action induces an action of $GL(t)$ on the coordinate ring on $M(m,t)\times M(t,n)$, which is a polynomial ring in $mt+tn$ variables. The classcial invariant theory studies the ideal of polynomials, which are invariant under this action. Let $\mu$ be the natural morphism of affine varieties $M(m,t)\times M(t,n)\lora M(m,n),$ $(A,B)\loma AB$, inducing a morphism of algebras
\bbs \mu^*:\O(M(m,n)\lora \O(M(m,t)\times M(t,n))\cong \O(M(m,t))\otimes \O(M(t,n)).\ees
Let $m^i_k$ (resp. ${n}^i_j, {p}^j_k$) be the standard generators of $\O(M(m,n))$, (resp. $\O(M(m,t)), \O(M(t,n))$), such that $\mu^*$ is given by
\bbs \mu^*(m^i_k)=\sum_k {n}^i_j\ot {p}^j_k.\ees
The first and the second fundamental theorems for general linear groups state that
\begin{enumerate}
\item Any invariant polynomial on $M(m,t)\times M(t,n)$ can be obtained by composing a polynimial on $M(m,n)$ with $\mu$, or, equivalently, the ideal of invariant polynomials is precisely $\Im\mu^*$. Thus, it is the quotient of $\O(M(m,n))$ by $\Ker\mu^*$.
\item The ideal $\Ker\mu^*$ in $\O(M(m,n))$ is generated by the minors of rank $(t+1)\times (t+1)$, i.e., it is the ideal $I_{t+1}$.\end{enumerate}

In this section we formulate and prove a quantum analogue of the above theorems for arbitrary Hecke operators $S,R$ and $T$, acting on $U,V$ and $W$ respectively. Thus $M(m,t)$ (resp. $M(t,n), M(m,n), GL(t)$) will be replaced by $\M_{RS}$ (resp. $\M_{TR},\M_{TS}, \Hb_R$). Here, $\Hb_R$ is a Hopf algebra associated to $R$. The left action in \rref{eqinv1} is replaced by a right coaction of $\Hb_R$ on $\M_{TR}\ot \M_{RS}$. 
The set of polynomials on $M(m,t)\times M(t,n)$, invariant with the action of $GL(t)$ now corresponds to the set of coinvariants of the coaction $\delta_{RST}$, i.e. the set of $x\in \M_{TR}\ot \M_{RS}$ such that $\delta_{RST}(x)=x\ot 1$.

In the quantum case, the are (at least) two ways to define the action of the Hopf algebra $\Hb_R$ on the algebra $\M_{RS}$, which coincide when the Hecke operators reduce to the ordinary flip operators. Therefore there are (at least) two versions of the fundamental theorems, depending on the way we define the coaction and on the way we define the algebra structure on $\M_{TR}\ot\M_{RS}$. In \ref{sect-firstversion} we define an algebra structure on $\M_{TR}\ot\M_{RS}$ in a usual way and choose an appropriate coaction of $\Hb_R$ on it. In this setting $\M_{TR}$ is an $\Hb_R$-comodule algebra, $\M_{RS}$ is an $\Hb_R$-comodule but not comodule algebra. Another setting is considered Subsection \ref{sect-secondversion}, in which we modify $\M_{RS}$ so that it becomes $\Hb_R$-comodule algebra (in fact, $R$ is replaced by $\hat R$). The algebra structure on $\M_{TR}\ot \M_{\hat RS}$ is also modified making the new algebra an $\Hb_R$-comodule algebra.

\subsection{Fundamental theorems for quatum groups of type $A$, the first version}\label{sect-firstversion}
In formulating the fundamental theorems for quantum groups of type $A$, we need to intoduce the multiplication map $\mu^*$ and the coaction of the quantum group $\Hb_R$.

We first define the morphism $\mu^*$. The linear map
\bbs \theta_1=\id\ot\db_V\ot \id :W^*\ot U\lora W^*\ot V\ot V^*\ot U\lora \T(W^*\ot V)\ot \T(V^*\ot U)\ees
 induces an algebra morphism
\bbs \theta:\T(W^*\ot U)\lora \T(W^*\ot V)\ot \T(V^*\ot U).\ees
Note that $\T(W^*\ot V)\ot \T(V^*\ot U)$ is the tensor product of the algebras$\T(W^*\ot V)$ and $ \T(V^*\ot U)$, that is, the elements from the latter algebras commute in $\T(W^*\ot V)\ot \T(V^*\ot U)$. In other words, we identify the two vector spaces
\bba\label{eq41} (W^*\ot V\ot V^*\ot U)^{\ot n}\cong (W^*\ot V)^{\ot n}\ot (V^*\ot U)^{\ot n},\eea
by means of the usual flip operator (that changes orders of tensor components). The restriction of $\theta$ on the $n^{\rm th}$ component is then the $n^{\rm th}$ tensor power of $\theta_1$, taking in account the above identification.

Fix bases of $U,V,W$ and then define their dual bases on $U^*,V^*,W^*$. These bases  define  bases for $W^*\ot U,W^*\ot V, V^*\ot U$, which will be denoted by $M=(m^i_k), N=({n}^i_j), L=(l^j_k),$ respectively. For convennience, we shall omit all tensor signs when describe an element of the algebra $\T(W^*\ot V)\ot \T(V^*\ot U)$. Then we have $\theta(M)=N\cdotot L$. Since in the algebra $\T(W^*\ot V)\ot \T(V^*\ot U)$, the entries of $N$ and $L$ commute, i.e. $L_1\cdot N_2=N_2\cdot L_1,$ where $L_1=L\ot \id, L_2=\id\ot L$, we have
\bbs \theta(M_1M_2)=(N\cdotot L)_1(N\cdotot L)_2=N_1{N}_2\cdotot L_1 L_2.\ees
Notice that $(M\cdotot N)_1=M_1\cdotot N_1, (M\cdotot N)_2=M_2\cdotot N_2$.
 
Combining $\theta$ with the quotient map $\T(W^*\ot V)\ot \T(V^*\ot U)\lora \M_{TR}\ot \M_{RS}$, we obtain an algebra morphism
\bbs \bar\theta:\T(W^*\ot U)\lora \M_{TR}\ot \M_{RS}.\ees
On $\M_{TR}\ot\M_{RS}$, we have
\bbs T(N\cdotot L)_1 (N\cdotot L)_2= TN_1N_2\cdotot L_1L_2
= N_1N_2R\cdotot L_1L_2
= N_1N_2\cdotot L_1L_2S
=(N\cdotot L)_1(N\cdotot L)_2S.\ees
Thus, $\bar\theta(TM_1M_2-M_1M_2S)=0$. Hence, it factorizes to a morphism
\bbs \mu^*:\M_{TS}\lora \M_{TR}\ot \M_{RS}.\ees

Next, we define the coaction of the Hopf algebra $\Hb_R$.
The Hopf algebra $\Hb_R$ is by definition the Hopf envelope of the bialgebra $\E_R$ \cite{manin1}, that is, there exists uniquely a bialgebra morphism $i:\E_R\lora \Hb_R$ such that any bialgebra morphism $f:\E_R\lora H$ to a Hopf algebra $H$ factorizes uniquely through as a composition of $i$ and a Hopf algebra morphism $j:\Hb_R\lora H$, $f=j\cdot i$. The Hopf envelope of any bialgebras exists, hence we can define $\Hb_R$ for any Hecke operators $R$. If $R$ satisfies the Yang-Baxter equation (so for example, when $R$ is a Hecke operator) it is known that $\E_R$ is a coquasitriangular bialgebra \cite{l-t,kassel}. However, in general, the coquasitriangular structure on $\E_R$ cannot be extended on $\Hb_R$. 

We shall assume that $R$ is a Hecke symmetry which means that the operator $R^\#:=(\ev_V\ot \id_{V\ot V^*})(\id_{V^*}\ot R\ot \id_{V^*})(\id_{V^*\ot V}\ot \db_V):V^*\ot V\lora V\ot V^*$ is invertible, this condition provides the coquasitriangular structure on $\E_R$ can be extended on $\Hb_R$, the antipode on $\Hb_R$ is bijective \cite{schauen,haya96}, and the map $i$ is injective \cite[Thm.~2.3.5]{ph97b}.

By means of the injective bialgebra morphism $i$, we identify $\E_R$ with a subbialgebra of $\Hb_R$. Each (simple) $\E_R$-comodule becomes then (simple) $\Hb_R$-comodule. Further, since the antipode of $\Hb_R$ is an anti-homomorphism of (co)algebras, a left (resp. right) $\E_R$-comodule becomes a right (left) $\Hb_R$ by composing the coaction with the antipode. In particular, since $(V^{\ot n})^*$ is a left $\E_R$-comodule, it is a right $\Hb_R$-comodule. The coaction is explicitly given as follows. Let $\{x_i\}_{i=1}^d$ be a basis of $V$ and $\{\xi^i\}_{i=1}^d$ be the dual basis for $V^*$. Then $\{e^j_i:=\xi^j\ot x_i\}_{i,j=1}^d$ form a multiplicative matrix in $\E_R$ and the left coaction of $\E_R$ on $V^*$ is $\lam_{V^*}(\xi^i)=\sum_k e^i_k\ot \xi^k.$ The right coaction of $\Hb_R$ on $(V^{\ot n})^*$ is then
\bba\label{eqhraction} \delta(\xi^{i_1}\ot \cdots\ot \xi^{i_n})=\sum_{k_1,\ldots,k_n}\xi^{k_1}\ot \cdots\ot\xi^{k_n}\ot \SS(e^{j_1}_{k_1}\cdots e^{j_n}_{k_n}), \; \mbox{ $\SS$ denotes the antipode}.\eea
Since $\SS$ is injective, we also have 
\bbs \End^{\Hb_R}((V^{\ot n})^*)= {}^{\E_R}\End((V^{\ot n})^*).\ees
In particular, $(\Rt)_w$ are endomorphisms of right $\Hb_R$-comodules, for all $w\in \Ss_n$. Therefore, $\M_{RS}$ is a right $\Hb_R$-subcomodule of $\T(V^*\ot U)$. Notice that $\T(V^*\ot U)$ and hence $\M_{RS}$, is not an $\Hb_R$-comodule algebra. The reason is that the usual isomorphism 
$ V^{*\ot n}\ot V^{*\ot m}\lora V^{*\ot m+n}$
is not an $\Hb_R$-comodule morphism.

We are now ready to formulate a quantum analogue of the first and the second fundamental theorems for general linear groups.
\begin{thm}\label{thm-inv1}
Let $S,T$ be Hecke operators $R$ be a Hecke symmetry with the parameter $q$ not a root of unity of order greater that 1. Let $\delta_{RST}$ be the coaction of $\Hb_R$ on $\M_{TR}\ot \M_{RS}$, which is the tensor product of the coactions of $\Hb_R$ on $\M_{TR}$ and $\M_{RS}$ given above. Then:
\begin{enumerate}\item The set of coinvariants in $\M_{TR}\ot \M_{RS}$ with respect to the coaction $\delta_{RST}$ is precisely $\Im\mu^*$.
\item The kernel of $\mu^*$ in $\M_{TS}$ is the ideal $I(\langle((r+1)^{s+1})\rangle )$ where $(r,s)$ is the birank of $S$.\end{enumerate}\end{thm}
\proof 
From Section \ref{sec3}, we know that $\M_{TS}$ decomposes into a direct sum of $\E_T^{\rm cop}\ot \E_S$-comodules and $\M_{TR}\ot \M_{RS}$ decomposes into a direct sum of $\E_T^{\rm cop}\ot \E_R\ot\E_R^{\rm cop}\ot \E_S$-comodules:
\bba\bbar{rcl}  \M_{TS}&\cong& \bigoplus_\lam W^*_\lam\ot U_\lam,\\
 \M_{TR}\ot \M_{RS}&\cong& \bigoplus_{\lam\mu} W^*_\lam\ot V_\lam\ot V^*_\mu\ot U_\mu.\eear\label{eq42}\eea

On the other hand, $V_\lam\ot V^*_\mu$ is a comodule over $\Hb_R$ and thus $\M_{TR}\ot \M_{RS}$ can be considered as an $\E_T^{\rm cop}\ot \Hb_R\ot \E_S$-comodule. Consider the trivial coaction of $\Hb_R$ on $\M_{TS}$, i.e., consider $\M_{TS}$ as a direct sum of copies of $\bK$ which is $\Hb_R$-comodule by mean of the unit element $\delta(1_\bK)=1_\bK\ot 1_{\Hb_R}$. 

\begin{lem}\label{lem-inv1} With the coaction of $\Hb_R$ described above, $\mu^*$ is a morphism of $\E_T^{\rm cop}\ot \Hb_R\ot \E_S$-comodules, the restriction of $\mu^*$ on $W^*_\lam\ot U_\lam$ can be given by the map $\db_{V_\lam}$:
\bbs \mu^*|_{W^*_\lam\ot U_\lam}=\id\ot\db_{V_\lam}\ot \id: W^*_\lam\ot U_\lam\lora W^*_\lam\ot V_\lam\ot V^*_\lam\ot U_\lam.\ees\end{lem}

Assume that Lemma is true. Since the map $\db_{V_\lam}:\bK\lora V_\lam\ot V^*_\lam$  is injective unless $V_\lam=0$, the restriction of $\mu$ on $W^*_\lam\ot U_\lam$ in injective unless $V_\lam=0$. This implies that $\Ker\mu^*$ is the set
$\displaystyle \bigoplus_{\lam,V_\lam=0}W^*_\lam\ot U_\lam$, which is precisely the set $\displaystyle \bigoplus_{\lam\subset((r+1)^{s+1})}W^*_\lam\ot U_\lam=I(\langle((r+1)^{s+1})\rangle ),$ by Theorem \ref{lem-vanish}. The first claim of Theorem is proved.

Since $V_\lam$ is simple over $\Hb_R$, 
\bbs \Hom_{\Hb_R}(\bK,V_\lam\ot V^*_\lam)=\Hom_{\Hb_R}(V_\lam,V_\lam)=\bK.\ees
Therefore, the image of $\db_{V_\lam}$ is the subspace of $\Hb_R$-invariants in $V_\lam\ot V^*_\lam$. Consequenlty, the set of invariants in $W^*_\lam\ot V_\lam\ot V^*_\lam\ot U_\lam$ is the image of $W^*_\lam\ot U_\lam$. On the other hand, if $\lam\neq \mu$ then
\bbs \Hom_{\Hb_R}(\bK,V_\lam\ot V^*_\mu)=\Hom_{\Hb_R}(V_\lam,V_\mu)=0.\ees
Hence, for $\lam\neq \mu$, the comodule  $W^*_\lam\ot V_\lam\ot V^*_\mu\ot U_\mu$ does not contains non-zero invariants. Taking the direct sum for all $\lam$ we prove the second claim of Theorem.

\vskip1ex
\noindent{\it Proof of Lemma \ref{lem-inv1}.}
Recall that the map $\theta$ is given in terms of the inclusion $W^*\ot U\lora \T(W^*\ot V)\ot \T(V^*\ot U)$, which in its order is given by the map $\db_V$. The restriction of $\theta$ on $(W^*\ot U)^{\ot n}$ is
\bbs \theta_n:W^{*\ot n}\ot U^{\ot n}\cong(W^*\ot U)^{\ot n}\stackrel{\theta_1^{\ot n}}{\lora} (W^*\ot V\ot V^*\ot U)^{\ot n}\cong
W^{*\ot n}\ot \Vn\ot\Vns\ot U^{\ot n}.\ees
Hence
\bbs \theta_n=\id_{V^{*\ot n}}\ot \db_{V^{\ot n}}\ot \id_{U^{\ot n}}:W^{*\ot n}\ot U^{\ot n}\lora W^{*\ot n}\ot \Vn\ot V^{*\ot n}\ot U^{\ot n}.\ees
Thus, is easy to see that $\theta$ is a morphism of $\E_T^{\rm cop}\ot \E_S$-comodules. On the other hand, with respect to the coaction of $\Hb_R$ on $(V^{\ot n})^*$ given in \rref{eqhraction}, $\db_{V^{\ot n}}$ is an $\Hb_R$-comodules morphism. Therefore $\theta_n$ is a morphism of $\E_T^{\rm cop}\ot \Hb_R\ot \E_S$-comodules.

The restriction of $\theta_n$ on $W^*_\lam\ot U_\mu$ is then the map
\bbs W^*_\lam\ot U_\mu\lora W^*_\lam\ot \Vn\ot\Vns \ot U_\mu\cong 
\bigoplus_\gamma W^*_\lam\ot V_\gamma\ot V^*_\gamma\ot U_\mu.\ees 
 The map $\mu^*$ is obtained by passing to quotients. According to Theorem \ref{thm1}, we can identify $\M_{TS}$ with a subspace of $\T(W^*\ot U)$ such that the quotient map $\T(W^*\ot U)\lora \M_{TS}$ is given by the sum $\Psi$ of the projectors $\Psi_n$. Thus $\mu^*$ can be considered as the composition $(\Psi_{TR}\ot \Psi_{RS})\theta\Psi_{TS}$. Therefore, the restriction of $\mu^*$ on $W^*_\lam\ot U_\lam$ is the map
\bbs \id\ot\db_{V_\lam}\ot \id:W^*_\lam\ot U_\lam\lora W^*_\lam\ot V_\lam\ot V^*_\lam\ot U_\lam.\ees
Since the map $\db_{V_\lam}$ is a morphism of $\Hb_R$-comodules, the above map is a morphism of $\E_T^{\rm cop}\ot \Hb_R\ot \E_S$-comodules. Lemma \ref{lem-inv1} is therefore proved.\eee

\subsection{Fundamental theorems for quantum groups of type $A$, the second version}\label{sect-secondversion}
 Since $R$ is a Hecke symmetry, $\Hb_R$ is a coquasitriangular Hopf algebra, i.e., the category of $\Hb_R$-comodules is braided (see, e.g., \cite[Chapter XI]{kassel}). We can modify the struture above so that the morphism $\mu^*$ is an $\Hb_R$-comodule algebra morphism. To do it, we identify the two vector spaces in \rref{eq41} by means of an $\Hb_R$-comodule isomorphism $\tau_{V^*V}:V^*\ot V\lora V\ot V^*$, $\xi^i\ot x_j\loma x_k\xi^l {R^{-1}}^{ik}_{jl}$. More precisely, in defining an isomorphism from $(W^*\ot V\ot V^*\ot U)^{\ot n}$ to $(W^*\ot V)^{\ot n}\ot (V^*\ot U)^{\ot n}$, whenever we have to interchange $V^*$ and $V$ we shall use the comodule isomorphism $\omega_{V^*V}$ above. We therefore modify the algebra structure on $\T(W^*\ot V)\ot\T(V^*\ot U)$ replacing the commuting relation of $N$ and $P$ by the following relation
\bba\label{eq43} L_1N_2=N_2R^{-1}PL_1,\eea
where $P$ is the matrix of the usual flip operator $x\ot y\lora y\ot x$, with respect to the basis $x_1,x_2,\ldots, x_d$, $P^{ij}_{kl}=\delta^i_l\delta^j_k$. Let $\omega_n$ denote the isomorphism from $(V\ot V^*)^{\ot n}\lora V^{\ot n}\ot V^{*\ot n}$ obtained by using $\omega_{V^*V}$, so, for example, $\omega_1=\id$, $\omega_2=\id\ot \omega_{V^*V}\ot \id$. 

Since the algebra structure on $\T(W^*\ot V)\ot\T(V^*\ot U)$ is modified, the map $\theta$ should also be modified.
The restriction of $\theta$ on $(W^*\ot U)^{\ot n}$ is now
\bbas &\theta_n:(W^*\ot U)^{\ot n}\stackrel{\theta_1^{\ot n}}{\lora} (W^*\ot V\ot V^*\ot U)^{\ot n}\stackrel{{\omega_n }}{\lora}
W^{*\ot n}\ot \Vn\ot\Vns\ot U^{\ot n},&\eeas
If we identify $(W^*\ot U)^{\ot n}$ with $W^{*\ot n}\ot U^{\ot n}$, then we can consider $\theta_n$ as a morphism
\bbas& \theta_n=\id_{V^{*\ot n}}\ot\omega_n{\db_V}^{\ot n}\ot \id_{U^{\ot n}}:W^{*\ot n}\ot U^{\ot n}\lora W^{*\ot n}\ot \Vn\ot V^{*\ot n}\ot U^{\ot n}.&\eeas
Thus, if we define the coaction of $\Hb_R$ on $V^{*\ot n}$ in such a way that this comodule is isomorphic to $(V^*)^{\ot n}$ as $\Hb_R$-comodules then $\omega_n{\db_V}^{\ot n}$ is a comodule morphism and hence so is $\theta_n$. Explicitly, the coaction is given by
\bbs \delta(\xi^{i_1}\ot\xi^{i_2}\ot\cdots\ot\xi^{i_n})=\sum_{k_1,k_2,\ldots,k_n}\xi^{k_1}\ot\xi^{k_2}\ot\cdots\ot\xi^{k_n}\ot\SS(e^{i_1}_{k_1})\SS(e^{i_2}_{k_2})\cdots\SS(e^{i_n}_{k_n}).\ees

Note that with respect to these coactions, the operators $\hat R_w, w\in \Ss_n$, where $\hat R=PRP$, are comodules endomorphisms of $V^{*\ot n}$. Hence the quotient $\M_{\hat RS}$ of $\T(V^*\ot U)$ is a comodule over $\Hb_R$. Further, since the usual identification $(V^*)^{\ot m}\ot (V^*)^{\ot n}\lora (V^*)^{\ot m+n}$ is a morphism of $\Hb_R$-comodules, $\M_{\hat RS}$ is an $\Hb_R$-comodule algebra.

On the other hand,  we can check that $\theta$ factorizes to an algebra morphism 
\bbs\mu^*_m:\M_{ST}\lora \M_{SR}\ot_{\sf m}\M_{\hat RT}\ees 
where in $\M_{SR}\ot_{\sf m}\M_{\hat RT}$, $N$ and $L$ commute by the rule in \rref{eq43}. Indeed, on $\M_{SR}\ot_{\sf m}\M_{\hat RT}$ we have
\bbas T(N\cdotot L)_1 (N\cdotot L)_2&=& TN_1N_2R^{-1}P\cdotot L_1L_2
= N_1N_2RR^{-1}P\cdotot L_1L_2
= N_1N_2R^{-1}P\cdotot\hat R L_1L_2\\
&=& N_1N_2R^{-1}P\cdotot L_1L_2S
=(N\cdotot L)_1 (N\cdotot L)_2S.\eeas

Thus, we obtain a coaction of $\Hb_R$ on $\M_{SR}\ot_{\sf m}\M_{\hat RT}$, for which the morphism $\mu^*_{\sf m}$ is an $\E_T^{\rm cop}\ot \Hb_R\ot \E_S$-comodule morphism. 

Although we have modified the algebra structure on $\M_{SR}\ot_{\sf m}\M_{\hat RT}$, this does not, affect the decomposition \rref{eq42}. That is, we still have an isomorphism of $\E_T^{\rm cop}\ot \Hb_R\ot \E_S$-comodules
\bbs 
 \M_{TR}\ot_{\sf m} \M_{RS}\cong \bigoplus_{\lam\mu} W^*_\lam\ot V_\lam\ot V^*_\mu\ot U_\mu.\ees
Therefore an anlogue of Lemma \ref{lem-inv1} can be easily obtained, whence one gets an analogue of Theorem \ref{thm-inv1}
\begin{thm}\label{thm-inv2}
Let $S,T$ be Hecke operators $R$ be a Hecke symmetry. Then $\mu^*_{\sf m}$ is an $\Hb_R$-comodule algebra morphism:
\begin{enumerate}\item The set of coinvariants in $\M_{TR}\ot_{\sf m} \M_{RS}$ with respect to the coaction $\delta_{RTS}$ is precisely $\Im\mu^*_{\sf m}$.
\item The kernel of $\mu^*_{\sf m}$ in $\M_{TS}$ is the ideal $I(\langle((r+1)^{s+1})\rangle )$ where $(r,s)$ is the birank of $T$.\end{enumerate}\end{thm}
\proof It remains to check that $\M_{SR}\ot_{\sf m}\M_{\hat RT}$ is an $\Hb_R$-comodule algebra. In fact, what we have done above is to define an algebra struture on the tensor product of two $\Hb_R$-comodule algebras $\M_{TR}$ and $\M_{\hat RS}$ (see, e.g., \cite{majid1}). Namely, the isomorphism $\omega_{V^*,V}$ gives rise to an $\Hb_R$-comodule isomorphism $\M_{\hat RS}\ot \M_{TR}\lora \M_{TR}\ot \M_{\hat RS}$. It follows from the standard argument that $\M_{TR}\ot_{\sf m} \M_{\hat RS}$ is an $\Hb_R$-comodule algebra.\eee

\section{Example: Standard quantum general linear groups}\label{sec7} The notion of quantum determinantal ideals presented here does not seem to have relationship with any quantum determinants. In fact, for an arbitrary Hecke operator, we aren't able to define any quantum determinant. However, in the case of standard $R$-matrix (see Subsection \ref{sec2.2}), the quantum minors are definable, and our notion of quantum determinantal ideals can be given in terms of these quantum minors. The quantum determinantal ideals associated to standard quantum matrix of type $A$ were studied by Goodearl, Lenagan and Rigal in \cite{gl,glr}, where the primeness was particularly proved. Here we show that in the case of standard $R$-matrix, our notion of quantum determinant ideal coincides with the notion give by Goodearl and Lenagan.

Recall that the standard $R$-matrix of type $A_{n-1}$ has, with respect to a certain basis $x_1,x_2,\ldots,x_n$, the following form:
\bbas &&{R_n}^{kl}_{ij}=
\frac{q^2-q^{2\va_{ij}}}{1+q^{2\va_{ij}}}\delta_{ij}^{kl}+\frac{q^{\va_{ij}}(q^2+1)}{1+q^{2\va_{ij}}}\delta_{ij}^{lk},\quad 1\leq i,j,k,l\leq n,\quad \va_{ij}:=\mbox{sign }(j-i)\eeas
The Hecke equation for $R_n$ is $(R_n-q^2)(R_n+1)=0$.
The quantum exterior algebra $\Lambda_{R_n}$ is the factor algebra of the non-commutative algebra $\bK\langle x_1,x_2,\ldots,x_n\rangle $ by the relations $x_ix_j=-qx_jx_i$ for $i\leq j$, it can be realized as subalgebra of $\bK\langle x_1,x_2,\ldots,x_n\rangle $ spanned by $q$-symmetrized tensors
$$x_{i_1}\wedge x_{i_2}\wedge\cdots\wedge x_{i_k}:=\sum_{\si\in\Ss_k}(-q)^{l(\si)} x_{i_1\si}\ot  x_{i_2\si}\ot\cdots\ot   x_{i_k\si},$$
for any sequence $(i_1< i_2<\ldots< i_k)$ of elements from $\{1,2,\ldots,n\}$, for $k=1,2,\ldots,n.$ In particular, ${\Lambda_{R_n}}_k$ is spanned by $x_{i_1}\wedge x_{i_2}\wedge\cdots\wedge x_{i_k}$.

Analogously, assume that ${\Lambda_{R_m}}_k$ has a basis consisting of $q$-antisymmetrized tensors $y_{j_1}\wedge y_{j_2}\wedge\cdots\wedge y_{j_k}$, for any sequence $(j_1< j_2<\ldots< j_k)$ of elements from $\{1,2,\ldots,m\}$, for $k=1,2,\ldots,m$. Then the space $({\Lambda_{R_m}}_k)^*=\Im\tau_k(E_{(1^k)})$ is canonically spanned by the set
$$\xi^{i_1}\wedge \xi^{i_2}\wedge\cdots\wedge \xi^{i_k}:=\sum_{\si\in\Ss_k}(-q)^{-l(\si)} \xi^{i_1\si}\ot  \xi^{i_2\si}\ot\cdots\ot   \xi^{i_k\si},$$
where $\xi^1,\xi^2,\ldots,xi^m$ is the dual basis to $y_1,y_2,\ldots, y_m$.

Denote $e_i^j:=\xi^j\ot x_i$. Then the algebra $\M_{R_mR_n}$ can be considered as a subspace of $\bK\langle e_1^1,\ldots,e_n^m\rangle $ spanned by the elements
\bbas e^{j_1j_2\ldots j_k}_{i_1i_2\ldots i_k}&:=&\sum_{\si,\tau\in\Ss_k}\frac{(-q)^{l(\si)}}{(-q)^{l(\si)}} e^{j_1\tau}_{i_1\si}\ot  e^{j_2\tau}_{i_2\si}\ot\cdots\ot   e^{j_k\tau}_{i_k\si}\\
&=&k! \sum_{\si\in\Ss_k}(-q)^{l(\si)} e^{j_1}_{i_1\si}\ot  e^{j_2}_{i_2\si}\ot\cdots\ot   e^{j_k}_{i_k\si}\eeas

The element $\sum_{\si\in\Ss_k}(-q)^{l(\si)} e^{j_1}_{i_1\si}\ot  e^{j_2}_{i_2\si}\ot\cdots\ot   e^{j_k}_{i_k\si}$ is precisely the quantum determinant of the submatrix of $Z=(e_i^j)$ formed on rows $i_1,i_2,\ldots, i_k$ and columns $j_1,j_2,\ldots, j_k$.

The fundamental theorems for standard quantum groups were proved by Goodearl, Lenagan and Rigal in \cite[Theorem 2.5]{gl} and \cite[Theorem 4.5]{glr}. We would like to mention that the setting of these theorems is slightly different from our setting here, namely, in defining the coaction of the quantum group and the algebra structure. In the language of our paper, the algebra structure of $\M_{TR}\ot \M_{RS}$ considered in [loc.cit] is the ordinary algebra structure, i.e., as in the setting of Subsection \ref{sect-firstversion}. The coaction of $\Hb_R$ on $\M_{RS}$ considered in [loc.cit] corresponds however to the coaction given in Subsection \ref{sect-secondversion}. One can do that because for standard deformation we have $\hat R_n:=PR_nP=R_n$.

\section{Example: The standard quantum general linear supergroups}\label{sec8}
There are two main differences between the coalgebras $\O(GL(n))$ and $\O(GL(m|n))$. Firstly, as a coalgebra, $\O(GL(n))$ is cosemisimple while $\O(GL(m|n))$ is not. Secondly, the determinant is a polynomial function while the super determinant is not. Looking at the more general construction of quantum groups of type $A$, we see that the quantum determinant is determined in terms of the quantum exterior algebra, which should be finite dimensional. Thus, if a Hecke opertor $R$ produces a finite dimensional quantum exterior algebra, we call the corresponding bialgebra $\E_R$ a quantum semi-group of type $A_n$, where $n$ is the of the Poincar\'e series of the quantum exterior algebra. In this case, $R$ is called even Hecke operator. It may happen that the quantum symmetric algebra has finite dimension, in this case $R$ is called odd Hecke operator \cite{gur1,ph97b}. A typical of non-even and non-odd Hecke operator is the flip operators in super geometry. Hence it is natural to suggest non-even non-odd Hecke operators define analogies of the general linear supergroups. It this then also natural to consider them as operators in the category of vector superspaces. It turns out, however, that the basis category does not play a great role. That is, it doesn't matter whether or not we consider $R$ as an operator in the category of vector superspaces and define $\E_R$ as a superbialgebra, many properties of $\E_R$ remain unchanged. In other words, many properties of $\E_R$ depend only on the intrinsic properties of $R$.

For example, let $R_d$ be the standard matrix considered in the previous section. If we assume that some of the basis vectors $x_1,x_2,\ldots,x_d$ have odd parity and the other have even parity, thus, $V$ is a vector superspace, then $R$ is an operator in the category of vector superspaces. In this category, the associated bialgebra $\E_R$ is defined differently but it remains cosemisimple. On the other hand, if we take a flip operator on a strict vector superspace and consider it as an operator in the category of (non-super) vector spaces, the associated bialgebra $\E_R$ remains non-cosemisimple.

All results in this paper hold in the category of vector superspaces. In fact, all what we have to do in the category of vector superspaces is to replace the ordinary flip operator by it super counterpart, i.e. to insert signs at some places.

Let now $R_{r|s}$ denote the super analogue of the standard $R$-matrices of type $A$. Explicitly, with respect to some homogeneous basis $x_1,x_2,\ldots,x_d$, $d=r+s$, where $\hat x_i:=\bar 0$ if $i\leq r$ and $\bar1$ if $i\geq r+1$, the operator $R_{r|s}$ has the following form:
\bbas &&{R_{r|s}}^{kl}_{ij}=
\frac{q^2-q^{2\va_{ij}}}{1+q^{2\va_{ij}}}\delta_{ij}^{kl}+\hat i\hat j\frac{q^{\va_{ij}}(q^2+1)}{1+q^{2\va_{ij}}}\delta_{ij}^{lk},\; 1\leq i,j,k,l\leq n=r+s,\; \va_{ij}:=\mbox{sign }(j-i),
\eeas
where $\hat i$ denotes the parity of $x_i$.
Then $R_{r|s}$ is a Hecke symmetry of birank $(r,s)$. The associated super bialgebra $\E_R$ and Hopf algebra $\Hb_R$ are called the function aglebras on the standard quantum super semi-group $M_q(r|s)$ and the standard quantum supergroup $GL_q(r|s)$, respectively, see, e.g. \cite{manin2,ph97d}.

Theorem \ref{thm-inv1} applied to this case gives us the fundamental theorems for standard quantum supergroups.
\begin{thm}\label{thm-inv3}
Let $M_q(m|n,r|s)$ denote the super bialgebra $\M_{R_{r|n}R_{r|s}}$ and let $GL_q(m|n)$ denote the Hopf superaglebra $\Hb_{R_{m|n}}$. Let $\mu^*$ be the algebra morphism $M_q(m|n,u|v)\lora M_q(m|n,r|s)\otimes M_q(r|s,u|v)$ induced from the map
\bbs e_j^i\lora \sum_{k=1}^{r+s}e^i_k\ot e^k_j\ees
where $\{e^i_j\}_{i=1,j=1}^{m+n\ u+v}$ is the standard generators of $M_q(m|n,u|v)$, similary, ${n}_k^i$ and ${p}^k_j$ are generators for $M_q(m|n,r|s)$ and $ M_q(r|s,u|v)$. Assume that $q$ is not a root ou unity of order greater than 1. Then, 
\begin{enumerate}\item the set of coinvariants of $M_q(m|n,r|s)\otimes M_q(r|s,u|v)$ with respect ot the coaction of $GL_q(r|s)$ is precisely $\Im\mu^*$,
\item the kernel of $\mu^*$ in $M_q(m|n,u|v)$ is the ideal $I(\langle((s+1)^{r+1})\rangle)=\sum_{\si_{u+1}\geq v+1}I_\si$.\end{enumerate}\end{thm}

Setting $q=1$ in the above theorem, we obtain the fundamental theorems for general linear super groups. These theorems can be formulated in the classical way. Let $\mu$ denote the multiplication of supermatrices
\bbs \mu:M(m|n,r|s)\times M(r|s,u|v)\lora M(m|n,r|s).\ees
Then $\mu$ induces a morphism $\mu^*$
\bbs \O(M(m|n,r|s))\lora \O(M(m|n,r|s)\times M(r|s,u|v))\cong \O(M(m|n,r|s))\otimes\O( M(r|s,u|v)).\ees
Consider the natural action of $GL(r|s)$ on $M(m|n,r|s)\times M(r|s,u|v)$: $(A,B)\loma (Ag^{-1},gB)$, which induces a natural coaction of $GL(r|s)$ on the function algebra on $M(m|n,r|s)\times M(r|s,u|v)$.
\begin{thm}\label{thm-inv4} We have:
\begin{enumerate}
\item A polynomial in $\O(M(m|n,r|s)\times M(r|s,u|v))$, invariant with the action of $GL(r|s)$, can be obtained by composing a polynomial in $\O(M(m|n,u|v))$ with $\mu$.
\item The kernel of $\mu^*$ is the ideal $I(\langle((s+1)^{r+1})\rangle)$.\end{enumerate}\end{thm}
\noindent{\it ~Remark.} Except for the case $n=s=v=0$, the ideal $I(\langle((s+1)^{r+1})\rangle)$ is not a determinantal ideal. It is an interesting problem to describe this ideal more explicitly.

\begin{center}\bf Acknowledgment\end{center}

\mythanks

\end{document}